\documentclass[12pt]{article}
\usepackage{amssymb}
\usepackage{amsmath}
\usepackage{makeidx}

\begin{document}

\renewcommand{\theequation}{\arabic{section}.\arabic{equation}}

\begin{center}
{\Huge Nonlinear parabolic flows with dynamic flux on the boundary}
\end{center}

\bigskip

\begin{center}
Viorel Barbu$^{a}$, Angelo Favini$^{\mbox{b}}$, Gabriela Marinoschi$^{\mbox{%
c,}\ast }$\footnotetext{%
Corresponding author
\par
E-mail addresses: vbarbu41@gmail.com (V. Barbu), favini@dm.unibo.it (A.
Favini), gabriela.marinoschi@acad.ro (G. Marinoschi)}
\end{center}

\bigskip

$^{a}$ {\small Al.I. Cuza University and Octav Mayer Institute of
Mathematics (Romanian Academy), Blvd. Carol I, No. 8, 700506 Ia\c{s}i,
Romania, e-mail: vbarbu41@gmail.com}

$^{b}${\small \ Department of Mathematics, University of Bologna, Piazza di
Porta S. Donato 5, 40126 Bologna, Italy, e-mail: favini@dm.unibo.it}

$^{c}$ {\small Institute of Mathematical Statistics and Applied Mathematics
of the Romanian Academy, Calea 13 Septembrie 13, 050711 Bucharest, Romania,
e-mail: gabriela.marinoschi@acad.ro}

\bigskip

{\small Abstract. A nonlinear divergence parabolic equation with dynamic
boundary conditions of Wentzell type is studied. The existence and
uniqueness of a strong solution is obtained as the limit of a finite
difference scheme, in the time dependent case and via a semigroup approach
in the time-invariant case.}

\medskip

{\small MSC 2010: 34A60, 35K55, 35K57, 35K61, 65N06}

{\small \medskip }

{\small Keywords: Nonlinear parabolic equations, dynamic boundary
conditions, reaction-diffusion, finite difference scheme, maximal monotone
operators}

\section{Introduction}

\setcounter{equation}{0}

This paper deals with the well-posedness of a nonlinear parabolic equation
posed in a bounded regular domain $\Omega $ (of class $C^{2},$ for instance)
of $\mathbb{R}^{N},$ $N\geq 1,$ coupled with a dynamic boundary condition of
reaction--diffusion type. More exactly, we study the problem
\begin{equation}
y_{t}-\nabla \cdot \beta (t,x,\nabla y)\ni f,\mbox{\ in }Q:=(0,T)\times
\Omega ,  \label{w1}
\end{equation}%
\begin{equation}
\beta (t,x,\nabla y)\cdot \nu +y_{t}\ni g,\mbox{\ on }\Sigma :=(0,T)\times
\Gamma ,  \label{w2}
\end{equation}%
\begin{equation}
y(0)=y_{0},\mbox{\ in }\Omega ,  \label{w3}
\end{equation}%
where $t\in (0,T),$ $T<\infty ,$ $x\in \Omega ,$ $\Gamma $ is the boundary
of $\Omega ,$ $\nu $ is the outward normal to $\Gamma ,$ $y_{t}=\frac{%
\partial y}{\partial t}$ and $\nabla y=\left\{ \frac{\partial y}{\partial
x_{i}}\right\} _{i=1,...,N}.$

In relation with various cases studied in this paper, certain combinations
of the following hypotheses will be used:

(H$_{1})$ For each $(t,x)\in \overline{Q},$ $\beta :\overline{Q}\times
\mathbb{R}^{N}\rightarrow \mathbb{R}^{N}$ is a maximal monotone graph with
respect to $r$ on $\mathbb{R}^{N}\times \mathbb{R}^{N},$ and it is derived
from a potential $j(t,x,r).$ The function $j$ is continuous on $\overline{Q}%
\times \mathbb{R}^{N},$ and for each $(t,x)\in \overline{Q},$ it is convex
with respect to $r.$ We denote%
\begin{equation}
\partial j(t,x,r)=\beta (t,x,r),\mbox{ for any }r\in \mathbb{R}^{N},\mbox{ }%
t\in \lbrack 0,T],\mbox{ }x\in \overline{\Omega },  \label{w4}
\end{equation}%
where $\partial j(t,x,\cdot )$ denotes the subdifferential of $j(t,x,\cdot
), $ that is
\begin{equation}
\partial j(t,x,r)=\{w\in \mathbb{R}^{N};\mbox{ }j(t,x,r)-j(t,x,\overline{r}%
)\leq w(r-\overline{r}),\mbox{ }\forall \overline{r}\in \mathbb{R}^{N}\}.
\label{w4-1}
\end{equation}%
Moreover, there is
\begin{equation}
\xi _{0}\in C(\overline{Q};\mathbb{R}^{N})\mbox{ with }\nabla \cdot \xi
_{0}\in L^{2}(\Omega ),\mbox{ such that }\xi _{0}(t,x)\in \beta (t,x,0),%
\mbox{ }\forall (t,x)\in \overline{Q}.  \label{w4-2}
\end{equation}

(H$_{2})$ (strong coercivity hypothesis): there exist $C_{i},C_{i}^{0}\in
\mathbb{R},$ $C_{1},C_{2}>0$, such that%
\begin{equation}
C_{1}\left\vert r\right\vert _{N}^{p}+C_{1}^{0}\leq j(t,x,r)\leq
C_{2}\left\vert r\right\vert _{N}^{p}+C_{2}^{0},\mbox{ }\forall (t,x)\in
\overline{Q},\mbox{ for }1<p<\infty .\mbox{ }  \label{w5}
\end{equation}%
(All over, $\left\vert \cdot \right\vert _{N}$ will denote the Euclidian
norm in $\mathbb{R}^{N}.)$

(H$_{3})$ (weak coercivity hypothesis) The functions $j$ and $j^{\ast }$
satisfy%
\begin{equation}
\lim_{\left\vert r\right\vert _{N}\rightarrow \infty }\frac{j(t,x,r)}{%
\left\vert r\right\vert _{N}}=+\infty ,\mbox{ uniformly with respect to
\thinspace }t,x,\mbox{\ }  \label{w293}
\end{equation}%
\begin{equation}
\lim_{\left\vert \omega \right\vert _{N}\rightarrow \infty }\frac{j^{\ast
}(t,x,\omega )}{\left\vert \omega \right\vert _{N}}=+\infty ,\mbox{
uniformly with respect to }t,x,  \label{w294}
\end{equation}%
where $j^{\ast }:\overline{Q}\times \mathbb{R}^{N}\rightarrow \mathbb{R}$ is
the conjugate of $j,$ defined by%
\begin{equation}
j^{\ast }(t,x,\omega )=\sup_{r\in \mathbb{R}^{N}}(\omega \cdot r-j(t,x,r)),%
\mbox{ for all }\omega \in \mathbb{R}^{N},\mbox{ }\forall (t,x)\in \overline{%
Q}.  \label{w290}
\end{equation}%
We note that (\ref{w293}) and (\ref{w294}) are equivalent with
\begin{equation}
\sup \left\{ \left\vert r\right\vert _{N};\mbox{ }r\in \beta
^{-1}(t,x,\omega ),\mbox{ }\left\vert \omega \right\vert _{N}\leq M\right\}
\leq W_{M},  \label{w293-0}
\end{equation}%
\begin{equation}
\sup \left\{ \left\vert \omega \right\vert _{N};\mbox{ }\omega \in \beta
(t,x,r),\mbox{ }\left\vert r\right\vert _{N}\leq M\right\} \leq Y_{M},
\label{w294-0}
\end{equation}%
respectively, where $M,$ $W_{M},$ $Y_{M}$ are positive constants.

(H$_{4})$ (symmetry at infinity) There exist $\gamma _{1},\gamma _{2}\geq 0$
such that
\begin{equation}
j(t,x,r)\leq \gamma _{1}j(t,x,-r)+\gamma _{2},\mbox{ }\gamma _{1}>0,\mbox{ }%
\gamma _{2}\geq 0.  \label{w296}
\end{equation}

(H$_{5})$ (regularity in $t)$ There exists $L>0$ such that
\begin{equation}
j(t,x,r)\leq j(s,x,r)+L\left\vert t-s\right\vert j(t,x,r),\mbox{ }\forall
t,s\in \lbrack 0,T],\mbox{ }x\in \overline{\Omega },\mbox{ }r\in \mathbb{R}%
^{N}.  \label{w5-1}
\end{equation}%
By (\ref{w4-2}) we see that subtracting $\xi _{0}\in \beta (t,x,0)$ from $%
\beta (t,x,r)$, and redefining $j(t,x,r)$ as $j(t,x,r)-j(t,x,0),$ we may
assume without loss of generality that
\begin{equation}
j(t,x,0)=0,\mbox{ }j(t,x,r)\geq 0,\mbox{ }j^{\ast }(t,x,r)\geq 0\mbox{ for
all }(t,x)\in \overline{Q},\mbox{ }r,\omega \in \mathbb{R}^{N}.  \label{w4-0}
\end{equation}

The strongly coercivity hypothesis (H$_{2}$) includes, for instance, the
situation
\begin{multline*}
\beta :\overline{Q}\times \mathbb{R}^{N}\rightarrow \mathbb{R}^{N},\mbox{ \ }%
\beta (t,x,r_{1},...,r_{N})=(\beta _{1}(t,x,r_{1}),...,\beta
_{N}(t,x,r_{N})), \\
\beta _{i}(t,x,r_{i})=\partial j_{i}(t,x,r_{i}),\mbox{ }i=1,...,N,
\end{multline*}%
where $j_{i}:\overline{Q}\times \mathbb{R}\rightarrow \mathbb{R}$ are
convex, continuous functions (with respect to $r_{i}),$ and continuous with
respect to $(t,x)\in \overline{Q},$ and $j:\overline{Q}\times \mathbb{R}%
^{N}\rightarrow \mathbb{R}$ is given by
\begin{equation*}
j(t,x,r_{1},...,r_{N})=j_{1}(t,x,r_{1})+...+j_{N}(t,x,r_{N}).
\end{equation*}%
For instance, one might take $j_{i}$ of the form
\begin{equation*}
j_{i}(t,x,r)=\alpha _{i}(t,x)\left\vert r\right\vert _{N}^{p}+\kappa
_{i}(t,x)\log (\left\vert r\right\vert _{N}+1)+\delta _{i}(t,x)\cdot
r+\delta _{i}^{1}(t,x),\mbox{ }i=1,...,N,\mbox{ }
\end{equation*}%
for $\alpha _{i},$ $\kappa _{i},$ $\delta _{i}^{1}\in C^{1}(\overline{Q}),$ $%
\kappa _{i}\geq 0,$ $\delta _{i}\in C^{1}(\overline{Q};\mathbb{R}^{N}).$

In particular, we get the parabolic equation with a non-isotropic $p$%
-Laplacian
\begin{eqnarray}
y_{t}-\sum_{i=1}^{N}\frac{\partial }{\partial x_{i}}\left( \alpha
_{i}(t,x)\left\vert \frac{\partial y}{\partial x_{i}}\right\vert ^{p-2}\frac{%
\partial y}{\partial x_{i}}\right) &=&f,\mbox{\ in }Q,  \notag \\
y_{t}+\sum_{i=1}^{N}\left( \alpha _{i}(t,x)\left\vert \frac{\partial y}{%
\partial x_{i}}\right\vert ^{p-2}\frac{\partial y}{\partial x_{i}}\right)
\cdot \nu _{i} &=&g,\mbox{\ on }\Sigma .  \label{w5-0}
\end{eqnarray}

More generally, we can consider instead of (\ref{w5-0}), a model for a
diffusion process in a fractured medium, described by the parabolic problem
\begin{eqnarray}
y_{t}-\sum_{i=1}^{N}\frac{\partial }{\partial x_{i}}\left( \left( \alpha
_{i}(t,x)+H\left( \frac{\partial y}{\partial x_{i}}-r_{i}\right) \right)
\left\vert \frac{\partial y}{\partial x_{i}}\right\vert ^{p-2}\frac{\partial
y}{\partial x_{i}}\right)  &\ni &f,\mbox{\ in }Q,  \notag \\
y_{t}+\sum_{i=1}^{N}\left( \left( \alpha _{i}(t,x)+H\left( \frac{\partial y}{%
\partial x_{i}}-r_{i}\right) \right) \left\vert \frac{\partial y}{\partial
x_{i}}\right\vert ^{p-2}\frac{\partial y}{\partial x_{i}}\right) \cdot \nu
_{i} &\ni &g,\mbox{\ on }\Sigma ,  \label{w5-00}
\end{eqnarray}%
where $r_{i}\in \mathbb{R},$ and $H$ is the Heaviside multivalued function, $%
H(s)=0$ for $s<0,$ $H(s)=[0,1]$ for $s=0,$ $H(s)=1$ for $s>1.$ In fact, a
discontinuous nondecreasing function $r\rightarrow \beta (t,x,r)$ becomes a
maximal monotone multivalued function by filling the jumps at the
discontinuity points $r_{i}$, that is by taking $\beta (t,x,r_{i})=[\beta
(t,x,r_{i}-0),\beta (t,x,r_{i}+0]$ and this is the natural way of treating
equation (\ref{w1}) with a discontinuous $\beta (t,x,\cdot ).$

Problem (\ref{w1})-(\ref{w3}) extends the classical Wentzell boundary
condition and models various phenomena in mathematical physics, and in
particular, diffusion and reaction--diffusion processes, phase-transition,
image restoring with observation on the boundary. If we view $%
E(y)=\int_{Q}j(t,x,\nabla y)dxdt$ as the energy of the system, then
hypothesis (H$_{2})$ describes diffusion processes with coercive and
differentiable energy, while (H$_{3})$ refers to systems with $W^{1,1}$
regular energy.

For various interpretations and treatment of the dynamic boundary conditions
(\ref{w2}), we refer e.g., to the works \cite{GRG-2006}, \cite{vazquez},
\cite{xiao}, \cite{coclite-favini-et-al}, \cite{favini-et-al}, \cite{frggr}.
In \cite{warma-2012} there are studied equations of the form $u_{t}-\nabla
\cdot (\left\vert \nabla u\right\vert ^{p-2}\nabla u)+\left\vert
u\right\vert ^{p-2}u+\alpha _{1}(x)=f$ in $Q,$ with the Wentzell boundary
condition $u_{t}-\nabla \cdot (\left\vert \nabla _{\Gamma }u\right\vert
^{p-2}\nabla _{\Gamma }u)+\left\vert \nabla u\right\vert ^{p-2}\partial
u/\partial \nu +\left\vert u\right\vert ^{p-2}u+\alpha _{2}(x,u)=g$ on $%
\Sigma .$ Previously, in \cite{gal-warma} and \cite{gal-2012-1} there were
studied problems with Wentzell boundary conditions of the form $u_{t}-\nabla
\cdot (a(\left\vert \nabla u\right\vert ^{2})\nabla u)+f(u)=h_{1}(x)$ (where
$a$ is a given nonnegative function), with the boundary condition $%
u_{t}+b(x)a(\left\vert \nabla u\right\vert ^{2})\partial u/\partial \nu
+g(u)=h_{2}(x).$

Compared with previous existence theory for problem (\ref{w1})-(\ref{w3}),
the novelty of the present work is two fold: the generality of the
nonlinearity $\beta ,$ which is discontinuous (that is, multivalued) and the
constructive approach based on a finite difference scheme, which permits to
treat the time dependent case.

The content of the paper is the following. In Section 3 we deal with the
strongly coercive case, under hypotheses (H$_{1}$), (H$_{2})$ and (H$_{5}).$
First, we prove the existence of a time-discretization solution to (\ref{w1}%
)-(\ref{w3}). Due to the generality assumed for $j$ we use a variational
principle involving an appropriate minimization problem. We also prove the
stability of the finite difference scheme. Then, we get the existence of a
weak solution as the strong limit of the $h$-discretized solution, with $h$
the time step. On the basis of some further arguments, it turns out that
this solution is strong and it is unique. In Section 4 we consider the
situation when $j$ is continuous and weakly coercive only, and exhibits a
symmetry at infinity, following hypotheses (H$_{1}$), (H$_{3})$-(H$_{5}).$
The latter case which provides a strong solution in the Sobolev space $%
W^{1,1}$ is in particular of interest in image processing with observation
on the boundary (see e.g. \cite{TB-VB-Coca}, \cite{vb-Rev-2013}). In Section
5 we present an alternative semigroup approach to the existence theory when
the potential $j$ is time independent$.$

\section{Notation and functional framework}

\setcounter{equation}{0}

In the following we denote by $\left\vert \cdot \right\vert _{N}$ the
Euclidian norm in $\mathbb{R}^{N},$ by $\left\vert \cdot \right\vert $ the
norm in $\mathbb{R},$ and by $u\cdot v$ the scalar product of $u,v\in
\mathbb{R}^{N}.$

Let $1\leq p\leq \infty .$ By $L^{p}(\Omega )$ we denote the space of $L^{p}$%
-Lebesgue integrable functions on $\Omega ,$ with the norm $\left\Vert \cdot
\right\Vert _{L^{p}(\Omega )}.$ Let $T>0.$ We set $Q:=(0,T)\times \Omega ,$ $%
\Sigma :=(0,T)\times \Gamma $ and denote by $L^{p}(Q)$ and $L^{p}(\Sigma )$
the corresponding $L^{p}$ spaces. We denote $W^{1,p}(\Omega )$ the Sobolev
spaces with the standard norm and $H^{1}(\Omega ):=W^{1,2}(\Omega ).$

We also denote by $p^{\prime }$ the conjugate of $p,$ that is $%
1/p+1/p^{\prime }=1.$

We define for $\sigma \in (0,1)$ and $p\in \lbrack 1,+\infty )$ the
fractional Sobolev spaces
\begin{equation*}
W^{\sigma ,p}(\Omega )=\left\{ z\in L^{p}(\Omega );\mbox{ }\frac{\left\vert
z(x)-z(x^{\prime })\right\vert }{\left\vert x-x^{\prime }\right\vert
^{\sigma +N/p}}\in L^{p}(\Omega \times \Omega )\right\}
\end{equation*}%
equipped with the natural norm (see, e.g., \cite{brezis-2011}, p. 314).

Next, for $s\in \mathbb{R},$ $s>1,$ not an integer, $s=m+\sigma ,$ $m$ being
the integer part of $s,$ one defines
\begin{equation*}
W^{s,p}(\Omega )=\left\{ z\in W^{m,p}(\Omega );\mbox{ }D^{\alpha }z\in
W^{\sigma ,p}(\Omega ),\mbox{ }\forall \alpha \mbox{ with }\left\vert \alpha
\right\vert =m\right\} .
\end{equation*}

If $z\in W^{1,r}(\Omega )$, with $r>1,$ it follows that the trace of $z$ on $%
\Gamma ,$ denoted by $\gamma (z),$ is well defined,
\begin{equation}
\gamma (z)\in W^{1-\frac{1}{r},r}(\Gamma ),\mbox{ \ }\left\Vert \gamma
(z)\right\Vert _{W^{1-\frac{1}{r},r}(\Gamma )}\leq C\left\Vert z\right\Vert
_{W^{1,r}(\Omega )},  \label{w21-4}
\end{equation}%
and the operator $z\rightarrow \gamma (z)$ is surjective from $%
W^{1,r}(\Omega )$ onto $W^{1-\frac{1}{r},r}(\Gamma )$ (see e.g., \cite%
{brezis-2011}, p. 315). We also have $\gamma (z)\in L^{1}(\Gamma )$ for $%
z\in W^{1,1}(\Omega ).$

For simplicity, when no confusion can be made, we still write $z$ instead of
$\gamma (z).$

Everywhere in the following, the gradient operator $\nabla ,$ as well as the
divergence $\nabla \cdot $ are considered in the sense of distributions on $%
\Omega .$

If $Y$ is a Banach space and $1\leq p\leq \infty ,$ we denote by $%
L^{p}(0,T;Y)$ the space of $L^{p}$ measurable $Y$-valued functions on $%
(0,T). $ By $W^{1,p}([0,T];Y)$ we denote the space $\left\{ y\in
L^{p}(0,T;Y);\frac{dy}{dt}\in L^{p}(0,T;Y)\right\} ,$ where $\frac{d}{dt}$
is considered in the sense of $Y$-valued distributions on $(0,T).$ Moreover,
each $y\in W^{1,p}([0,T];Y)$ is $Y$-valued absolutely continuous on $[0,T]$
and $\frac{d}{dt}$ exists a.e. on $(0,T)$ (see e.g. \cite{vb-springer-2010},
p. 23).

\section{The strongly coercive case}

\setcounter{equation}{0}

In this section we assume that hypotheses (H$_{1}),$ (H$_{2}),$ (H$_{5})$
hold, and $p>1.$ Let us define the space
\begin{equation*}
U=\{z\in L^{2}(\Omega );\mbox{ }\nabla z\in L^{p}(\Omega ),\mbox{ }\gamma
(z)\in L^{2}(\Gamma )\},
\end{equation*}%
endowed with the natural norm $\left\Vert z\right\Vert _{U}=\left\Vert
z\right\Vert _{L^{2}(\Omega )}+\left\Vert \nabla z\right\Vert _{L^{p}(\Omega
)}+\left\Vert \gamma (z)\right\Vert _{L^{2}(\Gamma )}.$

Recalling the Sobolev embeddings (\cite{brezis-2011}, p. 284), $%
W^{1,2}(\Omega )\subset L^{p^{\ast }}(\Omega ),$ if $N>2,$ where $p^{\ast }=%
\frac{2N}{N-2},$ $W^{1,2}(\Omega )\subset L^{q}(\Omega ),$ if $N=2,$ for any
$q\in \lbrack 2,+\infty ),$ $W^{1,2}(\Omega )\subset L^{\infty }(\Omega ),$
if $N=1,$ with continuous injections, we conclude that if $z\in U,$ we have
\begin{equation}
z\in W^{1,\overline{p}}(\Omega ),\mbox{ \ }\overline{p}=\left\{
\begin{array}{l}
p^{\ast },\mbox{ if }N>2,\mbox{ }p>p^{\ast }>2, \\
p,\mbox{ \ otherwise.}%
\end{array}%
\right.  \label{w21-5}
\end{equation}

In particular, if $p\geq 2,$ it follows that $U\subset H^{1}(\Omega )$ with
a dense and continuous embedding.

We mention for later use, that under assumption (\ref{w5}), one can easily
deduce that
\begin{equation}
\left\vert \xi \right\vert \leq C_{3}\left\vert r\right\vert
^{p-1}+C_{3}^{0},\mbox{ for any }\xi \in \beta (t,x,r),  \label{w21-2}
\end{equation}%
where $C_{3}$ and $C_{3}^{0}$ are positive constants.

We also assume that
\begin{equation}
y_{0}\in U,\mbox{ }f\in L^{2}(Q),\mbox{ }g\in L^{2}(\Sigma ).  \label{w21-3}
\end{equation}

\subsection{Existence and stability of the solution to the time-discretized
system}

We consider an equidistant partition $0=t_{0}\leq t_{1}\leq t_{2}\leq
...\leq t_{n}=T$ of the interval $[0,T],$ with $t_{i}=ih$ for $i=1,...,n,$ $%
h=T/n,$ and the finite sequences $\{f_{i}^{h}\}_{i=1,...,n},$ $%
\{g_{i}^{h}\}_{i=1,...,n},$ defined by the time averages
\begin{equation}
f_{i}^{h}=\frac{1}{h}\int_{(i-1)h}^{ih}f(s)ds,\mbox{ \ }g_{i}^{h}=\frac{1}{h}%
\int_{(i-1)h}^{ih}g(s)ds.  \label{w11}
\end{equation}%
We note that $f_{i}^{h}\in L^{2}(\Omega ),$ $g_{i}^{h}\in L^{2}(\Gamma ).$
We consider the time discretized system
\begin{equation}
\frac{y_{i+1}^{h}-y_{i}^{h}}{h}-\nabla \cdot \beta (t_{i+1},x,\nabla
y_{i+1}^{h})\ni f_{i+1}^{h},\mbox{\ in }\Omega ,\mbox{ }i=0,...,n-1,
\label{w12}
\end{equation}%
\begin{equation}
\beta (t_{i+1},x,\nabla y_{i+1}^{h})\cdot \nu +\frac{y_{i+1}^{h}-y_{i}^{h}}{h%
}\ni g_{i+1}^{h},\mbox{\ on }\Gamma ,  \label{w13}
\end{equation}%
\begin{equation}
y_{0}^{h}=y_{0},\mbox{\ in\ }\Omega .  \label{w14}
\end{equation}

\medskip

\noindent \textbf{Definition 2.1.} We call a \textit{weak solution} to the
time-discretized system (\ref{w12})-(\ref{w13}), a set of functions $%
\{y_{i}^{h}\}_{i=1,...,n},$ $y_{i}^{h}\in U,$ which satisfies (for each $%
i=1,...,n-1)$
\begin{eqnarray}
&&\int_{\Omega }y_{i+1}^{h}\psi dx+h\int_{\Omega }\eta _{i+1}^{h}\cdot
\nabla \psi dx+\int_{\Gamma }y_{i+1}^{h}\psi d\sigma  \label{w15} \\
&=&\int_{\Omega }y_{i}^{h}\psi dx+\int_{\Gamma }y_{i}^{h}\psi d\sigma
+h\int_{\Omega }f_{i+1}^{h}\psi dx+h\int_{\Gamma }g_{i+1}^{h}\psi d\sigma ,%
\mbox{ for any }\psi \in U,  \notag
\end{eqnarray}%
for some measurable function $\eta _{i+1}^{h}$, such that $\eta
_{i+1}^{h}(x)\in \beta (t_{i+1},x,\nabla y_{i+1}^{h}(x)),$ a.e. $x\in \Omega
.$

\medskip

We mention that in the third integral on the left-hand side in (\ref{w15})
we understand by $y_{i+1}^{h}$ the trace of $y_{i+1}^{h}\in U$ on $\Gamma $.
If $y_{i}^{h}$ is a solution, it follows by (\ref{w5}) and (\ref{w21-2})
that
\begin{equation}
j(t_{i},\cdot ,\nabla y_{i}^{h})\in L^{1}(\Omega ),\mbox{ }\eta _{i}^{h}\in
(L^{p^{\prime }}(\Omega ))^{N},\mbox{ }i=0,...,n.  \label{w14-0}
\end{equation}%
Proposition 2.2 is concerned with the stability of the finite difference
scheme (\ref{w12})-(\ref{w13}).

\medskip

\noindent \textbf{Proposition 2.2.} \textit{Let us assume} (\ref{w21-3})
\textit{and} $j(0,\cdot ,\nabla y_{0})\in L^{1}(\Omega ).$ \textit{System} (%
\ref{w12})-(\ref{w13}) \textit{has a unique weak solution satisfying}
\begin{equation}
\left\Vert y_{i}^{h}\right\Vert _{L^{2}(\Omega )}\leq C,\mbox{ \ }i=1,...,n,
\label{w17}
\end{equation}%
\begin{equation}
\left\Vert \gamma (y_{i}^{h})\right\Vert _{L^{2}(\Gamma )}\leq C,\mbox{ \ }%
i=1,...,n,  \label{w18}
\end{equation}%
\begin{equation}
h\sum_{i=0}^{m-1}\left\Vert \nabla y_{i+1}^{h}\right\Vert _{L^{p}(\Omega
)}^{p}\leq C,\mbox{ }m=1,...,n,  \label{w18-0}
\end{equation}%
\begin{equation}
h\sum_{i=0}^{m-1}\left\Vert \frac{y_{i+1}^{h}-y_{i}^{h}}{h}\right\Vert
_{L^{2}(\Omega )}^{2}\leq C,\mbox{ \ }m=1,...,n,  \label{w19}
\end{equation}%
\begin{equation}
h\sum_{i=0}^{m-1}\left\Vert \frac{\gamma (y_{i+1}^{h})-\gamma (y_{i}^{h})}{h}%
\right\Vert _{L^{2}(\Gamma )}^{2}\leq C,\mbox{ \ }m=1,...,n,  \label{w19-0}
\end{equation}%
\begin{equation}
h\sum_{i=0}^{m-1}\int_{\Omega }j(t_{i+1},x,\nabla y_{i+1}^{h})dx\leq C,\mbox{
\ }i=1,...,n,  \label{w20}
\end{equation}%
\textit{where }$C$\textit{\ is a positive constant, independent of }$h$%
\textit{. }

\medskip

\noindent \textbf{Proof. }Let us fix $t\in \lbrack 0,T],$ $w_{1}\in
L^{2}(\Omega ),$ $w_{2}\in L^{2}(\Gamma )$ and consider the intermediate
problem%
\begin{eqnarray}
u-h\nabla \cdot \beta (t,x,\nabla u) &\ni &w_{1}\mbox{ in }\Omega ,
\label{w20-0} \\
u+h\beta (t,x,\nabla u)\cdot \nu &\ni &w_{2}\mbox{ on }\Gamma .  \notag
\end{eqnarray}%
We define $b\in U^{\prime }$ (the dual of the space $U)$ by%
\begin{equation}
b(\psi ):=\int_{\Omega }w_{1}(x)\psi (x)dx+\int_{\Gamma }w_{2}(\sigma )\psi
(\sigma )d\sigma ,\mbox{ for all }\psi \in U,  \label{w20-00}
\end{equation}%
and note that
\begin{equation}
\left\vert b(\psi )\right\vert \leq \left\Vert w_{1}\right\Vert
_{L^{2}(\Omega )}\left\Vert \psi \right\Vert _{L^{2}(\Omega )}+\left\Vert
w_{2}\right\Vert _{L^{2}(\Gamma )}\left\Vert \psi \right\Vert _{L^{2}(\Gamma
)}\mbox{ for all }\psi \in U.  \label{w20-1}
\end{equation}%
For $u\in U$ it is clear by (\ref{w5}) and (\ref{w21-2}) that, for $t$
fixed,
\begin{equation*}
j(t,\cdot ,\nabla u)\in L^{1}(\Omega ),\mbox{ }\eta \in (L^{p^{\prime
}}(\Omega ))^{N},\mbox{ }
\end{equation*}%
for all measurable sections $\eta (x)$ of $\beta (t,x,\nabla u(x)).$

We call a weak solution to (\ref{w20-0}) a function $u\in U,$ such that
there is $\eta \in (L^{p^{\prime }}(\Omega ))^{N},$ $\eta (x)\in \beta
(t,x,\nabla u(x))$ a.e. $x\in \Omega ,$ and%
\begin{equation}
\int_{\Omega }(u\psi +h\eta \cdot \nabla \psi )dx+\int_{\Gamma }u\psi
d\sigma =b(\psi ),\mbox{ for all }\psi \in U.  \label{w22-1}
\end{equation}

To prove that (\ref{w20-0}) has a solution we use a variational argument,
i.e., we show that a solution to this equation is retrieved as a solution to
the minimization problem
\begin{equation}
\mbox{Min }\left\{ \varphi (u);\mbox{ }u\in U\right\} ,  \label{w26}
\end{equation}%
where $\varphi :L^{2}(\Omega )\rightarrow \mathbb{R\cup \{+\infty \}}$ is
given by
\begin{equation}
\varphi (u)=\left\{
\begin{array}{l}
\frac{1}{2}\int_{\Omega }u^{2}dx+h\int_{\Omega }j(t,x,\nabla u)dx+\frac{1}{2}%
\int_{\Gamma }u^{2}d\sigma -b(u),\mbox{ if }u\in U,\mbox{ } \\
+\infty ,\mbox{ \ \ \ \ \ \ \ \ \ \ \ \ \ \ \ \ \ \ \ \ \ \ \ \ \ \ \ \ \ \
\ \ \ \ \ \ \ \ \ \ \ \ \ \ \ \ \ \ \ \ \ \ \ \ \ \ \ \ otherwise.}%
\end{array}%
\right.  \label{w25}
\end{equation}%
By (\ref{w5}) and (\ref{w20-1}) we have
\begin{eqnarray}
\varphi (u) &\geq &\frac{1}{2}\left\Vert u\right\Vert _{L^{2}(\Omega
)}^{2}+hC_{1}\left\Vert \nabla u\right\Vert _{L^{p}(\Omega )}^{p}+hC_{1}^{0}+%
\frac{1}{2}\left\Vert u\right\Vert _{L^{2}(\Gamma )}^{2}-\left\vert
b(u)\right\vert  \notag \\
&\geq &\frac{1}{4}\left\Vert u\right\Vert _{L^{2}(\Omega
)}^{2}+hC_{1}\left\Vert \nabla u\right\Vert _{L^{p}(\Omega )}^{p}+\frac{1}{4}%
\left\Vert u\right\Vert _{L^{2}(\Gamma )}^{2}  \label{w25-0} \\
&&+hC_{1}^{0}-4\left\Vert w_{1}\right\Vert _{L^{2}(\Omega )}^{2}-4\left\Vert
w_{2}\right\Vert _{L^{2}(\Gamma )}^{2},\mbox{ }\forall u\in U.  \notag
\end{eqnarray}

It is also easily seen that $\varphi $ is proper, strictly convex and lower
semicontinuous (l.s.c. for short) on $L^{2}(\Omega ).$ Let us denote by $%
d=\inf\nolimits_{u\in U}\varphi (u)$ and let us consider a minimizing
sequence $\{u_{n}\}_{n\geq 1}$ for $\varphi .$ Then, we have
\begin{equation}
d\leq \varphi (u_{n})=\frac{1}{2}\int_{\Omega }u_{n}^{2}dx+h\int_{\Omega
}j(t,x,\nabla u_{n})dx+\frac{1}{2}\int_{\Gamma }u_{n}^{2}d\sigma
-b(u_{n})\leq d+\frac{1}{n}.  \label{w27}
\end{equation}%
By (\ref{w25-0}) it follows that
\begin{equation*}
\left\Vert u_{n}\right\Vert _{L^{2}(\Omega )}^{2}+h\left\Vert \nabla
u_{n}\right\Vert _{L^{p}(\Omega )}^{p}+\left\Vert u_{n}\right\Vert
_{L^{2}(\Gamma )}^{2}\leq C,\mbox{ }\forall n\in \mathbb{N},
\end{equation*}%
where $C$ is a positive constant independent of $n.$ Therefore we can select
a subsequence $(n\rightarrow \infty )$ such that
\begin{equation*}
u_{n}\rightarrow u\mbox{ weakly in }L^{2}(\Omega ),
\end{equation*}%
\begin{equation*}
\gamma (u_{n})\rightarrow \chi \mbox{ weakly in }L^{2}(\Gamma ),\mbox{ as }%
n\rightarrow \infty ,
\end{equation*}%
\begin{equation*}
\nabla u_{n}\rightarrow \xi \mbox{ weakly in }(L^{p}(\Omega ))^{N},\mbox{ as
}n\rightarrow \infty .
\end{equation*}%
It follows that $\xi =\nabla u,$ a.e. in $\Omega ,$ and so $u\in W^{1,%
\overline{p}}(\Omega )$ with $\overline{p}$ given by (\ref{w21-5}). This
implies that $\gamma (u)\in W^{1-\frac{1}{\overline{p}},\overline{p}}(\Omega
)$, and it is clear that $\chi =\gamma (u)$ a.e. on $\Gamma .$ Then, $u\in
U. $ Moreover, by (\ref{w21-2}) we have also $\eta \in (L^{p^{\prime
}}(\Omega ))^{N},$ where $\eta (x)\in \beta (t,x,\nabla u(x))$ a.e. $x\in
\Omega .$ Since $\varphi $ is convex and continuous it is also weakly l.s.c.
on $L^{2}(\Omega )$ and so $\liminf\limits_{n\rightarrow \infty }\varphi
(u_{n})\geq \varphi (u).$ Passing to the limit in (\ref{w27}), as $%
n\rightarrow \infty ,$ we get that $\varphi (u)=d,$ as claimed.

Next, we connect this solution to the solution to (\ref{w20-0}). Let $%
\lambda >0$ and define the variation $u^{\lambda }=u+\lambda \psi ,$ for all
$\psi \in C^{\infty }(\overline{Q}).$ We have $\varphi (u)\leq \varphi
(u^{\lambda }),$ for any $\lambda >0.$ Replacing the expression of $\varphi
, $ dividing by $\lambda $ and letting $\lambda \rightarrow 0,$ we get
\begin{equation*}
\int_{\Omega }(u\psi +h\eta \cdot \nabla \psi )dx+\int_{\Gamma }u\psi
d\sigma -b(\psi )\leq 0,
\end{equation*}%
for all $\psi \in C^{\infty }(\overline{Q}).$ By density, this extends to
all $\psi \in U$. Changing $\psi $ to $-\psi $ and making the same calculus
we obtain the converse inequality, so that in conclusion we find that the
solution to (\ref{w26}) satisfies (\ref{w22-1}).

Relying on this result we deduce in an iterative way that system (\ref{w12}%
)-(\ref{w13}) has a unique weak solution. We observe that it can be
rewritten as
\begin{equation}
\int_{\Omega }y_{i+1}^{h}\psi dx+h\int_{\Omega }\eta _{i+1}^{h}\cdot \nabla
\psi dx+\int_{\Gamma }y_{i+1}^{h}\psi d\sigma =b_{i+1}(\psi ),\mbox{ }%
\forall \psi \in U,  \label{w23}
\end{equation}%
for $i=0,...n-1,$ where $\eta _{i+1}^{h}(x)\in \beta (t_{i+1},x,\nabla
y_{i+1}^{h}(x))$ a.e. $x\in \Omega $ and $b_{i+1}\in U^{\prime }$ is given
by
\begin{equation}
b_{i+1}(\psi ):=\int_{\Omega }(y_{i}^{h}\psi +hf_{i+1}^{h}\psi
)dx+\int_{\Gamma }(y_{i}^{h}+hg_{i+1}^{h})\psi d\sigma ,\mbox{ }\forall \psi
\in U,  \label{w24}
\end{equation}%
and%
\begin{multline*}
\left\vert b_{i+1}(\psi )\right\vert \leq \left( \left\Vert
y_{i}^{h}\right\Vert _{L^{2}(\Omega )}+h\left\Vert f_{i+1}^{h}\right\Vert
_{L^{2}(\Omega )}\right) \left\Vert \psi \right\Vert _{L^{2}(\Omega )} \\
+\left( \left\Vert y_{i}^{h}\right\Vert _{L^{2}(\Gamma )}+h\left\Vert
g_{i+1}^{h}\right\Vert _{L^{2}(\Gamma )}\right) \left\Vert \psi \right\Vert
_{L^{2}(\Gamma )},\mbox{ for }i=0,...,n-1.
\end{multline*}%
We begin with the equation (\ref{w23}) for $i=0,$ in which $b_{1}(\psi )$
satisfies the previous relation for $i=0.$ Setting $b=b_{1}$ in $\varphi $
we get that the corresponding problem (\ref{w26}) has a unique weak solution
$y_{1}^{h}$ which verifies (\ref{w23}) for $i=0$. Next, we set $b=b_{i+1}$
in $\varphi $ and by recurrence, we obtain a sequence of solutions $%
y_{i+1}^{h}\in U,$ which satisfy (\ref{w23}) for all $i=1,...,n-1.$ In
particular, $y_{i}^{h}\in W^{1,\overline{p}}(\Omega )$, with $\overline{p}$
given by (\ref{w21-5}).

To obtain the first estimate (\ref{w17}) we set $\psi =y_{i+1}^{h}$ in (\ref%
{w15}) and use (\ref{w4}) getting
\begin{eqnarray}
&&\int_{\Omega }(y_{i+1}^{h})^{2}dx+h\int_{\Omega }(j(t_{i+1},x,\nabla
y_{i+1}^{h})-j(t_{i+1},x,0))dx+\int_{\Gamma }(y_{i+1}^{h})^{2}d\sigma  \notag
\\
&\leq &\int_{\Omega }(y_{i+1}^{h})^{2}dx+h\int_{\Omega }\eta _{i+1}^{h}\cdot
\nabla y_{i+1}^{h}dx+\int_{\Gamma }(y_{i+1}^{h})^{2}d\sigma  \label{w24-0} \\
&=&\int_{\Omega }y_{i}^{h}y_{i+1}^{h}dx+\int_{\Gamma
}y_{i}^{h}y_{i+1}^{h}d\sigma +h\int_{\Omega
}f_{i+1}^{h}y_{i+1}^{h}dx+h\int_{\Gamma }g_{i+1}^{h}y_{i+1}^{h}d\sigma .
\notag
\end{eqnarray}%
Then, we sum up (\ref{w24-0}) from $i=0$\ to $i=m-1\leq n-1$, use (\ref{w4-0}%
),%
\begin{eqnarray}
&&\frac{1}{2}\left\Vert y_{m}^{h}\right\Vert _{L^{2}(\Omega
)}^{2}+h\sum_{i=0}^{m-1}\int_{\Omega }j(t_{i+1},x,\nabla y_{i+1}^{h})dx+%
\frac{1}{2}\left\Vert y_{m}^{h}\right\Vert _{L^{2}(\Gamma )}^{2}
\label{w24-1} \\
&\leq &\frac{h}{2}\sum_{i=0}^{m-1}\left\Vert y_{i+1}^{h}\right\Vert
_{L^{2}(\Omega )}^{2}+\frac{h}{2}\sum_{i=0}^{m-1}\left\Vert
f_{i+1}^{h}\right\Vert _{L^{2}(\Omega )}^{2}+\frac{h}{2}\sum_{i=0}^{m-1}%
\left\Vert y_{i+1}^{h}\right\Vert _{L^{2}(\Gamma )}^{2}  \notag \\
&&+\frac{h}{2}\sum_{i=0}^{m-1}\left\Vert g_{i+1}^{h}\right\Vert
_{L^{2}(\Gamma )}^{2}+\frac{1}{2}\left\Vert y_{0}^{h}\right\Vert
_{L^{2}(\Omega )}^{2}+\frac{1}{2}\left\Vert y_{0}^{h}\right\Vert
_{L^{2}(\Gamma )}^{2},  \notag
\end{eqnarray}%
and obtain (since $j$ is continuous on $\overline{Q})$
\begin{eqnarray}
&&\left\Vert y_{m}^{h}\right\Vert _{L^{2}(\Omega )}^{2}+\left\Vert
y_{m}^{h}\right\Vert _{L^{2}(\Gamma )}^{2}+2C_{1}h\sum_{i=0}^{m-1}\left\Vert
\nabla y_{i+1}^{h}\right\Vert _{L^{p}(\Omega )}^{p}  \label{w28} \\
&\leq &C_{0}+h\sum_{i=1}^{m}\left( \left\Vert y_{i}^{h}\right\Vert
_{L^{2}(\Omega )}^{2}+\left\Vert y_{i}^{h}\right\Vert _{L^{2}(\Gamma
)}^{2}\right) ,  \notag
\end{eqnarray}%
where
\begin{equation*}
C_{0}=h\int_{0}^{T}(\left\Vert f(t)\right\Vert _{L^{2}(\Omega
)}^{2}+\left\Vert g(t)\right\Vert _{L^{2}(\Gamma )}^{2})dt+\left\Vert
y_{0}\right\Vert _{L^{2}(\Omega )}^{2}+\left\Vert y_{0}\right\Vert
_{L^{2}(\Gamma )}^{2}+2T\left\vert C_{1}^{0}\right\vert \mbox{meas}(\Omega ).
\end{equation*}%
Using a variant of the discrete Gronwall's lemma (see e.g., \cite{gm-cc-08})
we get
\begin{equation}
\left\Vert y_{m}^{h}\right\Vert _{L^{2}(\Omega )}^{2}+\left\Vert
y_{m}^{h}\right\Vert _{L^{2}(\Gamma )}^{2}\leq 2e^{T}(\left\Vert
y_{0}^{h}\right\Vert _{L^{2}(\Omega )}^{2}+\left\Vert y_{0}^{h}\right\Vert
_{L^{2}(\Gamma )}^{2}+C_{0}),  \label{w29}
\end{equation}%
\begin{equation}
h\sum_{i=1}^{m}\left\Vert y_{i}^{h}\right\Vert _{L^{2}(\Omega
)}^{2}+h\sum_{i=1}^{m}\left\Vert y_{i}^{h}\right\Vert _{L^{2}(\Gamma
)}^{2}\leq e^{T}(\left\Vert y_{0}\right\Vert _{L^{2}(\Omega
)}^{2}+\left\Vert y_{0}^{h}\right\Vert _{L^{2}(\Gamma )}^{2}+C_{0}).
\label{w30}
\end{equation}%
By the last two relations, (\ref{w28}) and (\ref{w24-1}) we obtain (\ref{w17}%
)-(\ref{w18-0}) and (\ref{w20}).

To get (\ref{w19}) we set $\psi =\frac{y_{i+1}^{h}-y_{i}^{h}}{h}$ in (\ref%
{w15}) and we obtain
\begin{eqnarray*}
&&\left\Vert \frac{y_{i+1}^{h}-y_{i}^{h}}{h}\right\Vert _{L^{2}(\Omega
)}^{2}+\left\Vert \frac{y_{i+1}^{h}-y_{i}^{h}}{h}\right\Vert _{L^{2}(\Gamma
)}^{2}+\int_{\Omega }\eta _{i+1}^{h}\cdot \nabla \frac{y_{i+1}^{h}-y_{i}^{h}%
}{h}dx \\
&=&\int_{\Omega }f_{i+1}^{h}\frac{y_{i+1}^{h}-y_{i}^{h}}{h}dx+\int_{\Gamma
}g_{i+1}^{h}\frac{y_{i+1}^{h}-y_{i}^{h}}{h}d\sigma \\
&\leq &\frac{1}{2}\left\Vert f_{i+1}^{h}\right\Vert _{L^{2}(\Omega )}^{2}+%
\frac{1}{2}\left\Vert \frac{y_{i+1}^{h}-y_{i}^{h}}{h}\right\Vert
_{L^{2}(\Omega )}^{2}+\frac{1}{2}\left\Vert g_{i+1}^{h}\right\Vert
_{L^{2}(\Gamma )}^{2}+\frac{1}{2}\left\Vert \frac{y_{i+1}^{h}-y_{i}^{h}}{h}%
\right\Vert _{L^{2}(\Gamma )}^{2}.
\end{eqnarray*}%
Further we use again (\ref{w4}) and sum up from $i=0$ to $i=m-1\leq n-1.$ We
obtain
\begin{eqnarray}
&&h\sum_{i=0}^{m-1}\left\Vert \frac{y_{i+1}^{h}-y_{i}^{h}}{h}\right\Vert
_{L^{2}(\Omega )}^{2}+h\sum_{i=0}^{m-1}\left\Vert \frac{y_{i+1}^{h}-y_{i}^{h}%
}{h}\right\Vert _{L^{2}(\Gamma )}^{2}\leq h\sum_{i=0}^{m-1}\left\Vert
f_{i+1}^{h}\right\Vert _{L^{2}(\Omega )}^{2}  \label{w30-1} \\
&&+h\sum_{i=0}^{m-1}\left\Vert g_{i+1}^{h}\right\Vert _{L^{2}(\Gamma
)}^{2}+2\sum_{i=0}^{m-1}\int_{\Omega }(j(t_{i+1},x,\nabla
y_{i}^{h})-j(t_{i+1},x,\nabla y_{i+1}^{h}))dx.  \notag
\end{eqnarray}%
By (\ref{w5-1}) we have
\begin{equation*}
j(t,x,r)\leq \frac{j(s,x,r)}{1-L\left\vert t-s\right\vert },\mbox{ for }%
t,s\in \lbrack 0,T],\mbox{ }\left\vert t-s\right\vert <1/L.
\end{equation*}%
Then, we compute
\begin{eqnarray*}
&&\sum_{i=0}^{m-1}\int_{\Omega }(j(t_{i+1},x,\nabla
y_{i}^{h})-j(t_{i+1},x,\nabla y_{i+1}^{h}))dx \\
&=&\sum_{i=0}^{m-1}\int_{\Omega }(j(t_{i+1},x,\nabla
y_{i}^{h})-j(t_{i},x,\nabla y_{i}^{h}))dx+\int_{\Omega }j(0,x,\nabla y_{0})dx
\\
&\leq &Lh\sum_{i=0}^{m-1}\int_{\Omega }j(t_{i+1},x,\nabla
y_{i}^{h})dx+\int_{\Omega }j(0,x,\nabla y_{0})dx \\
&\leq &\frac{Lh}{1-Lh}\sum_{i=0}^{m-1}\int_{\Omega }j(t_{i},x,\nabla
y_{i}^{h})dx+\int_{\Omega }j(0,x,\nabla y_{0})dx\leq C,\mbox{ }
\end{eqnarray*}%
for $h$ sufficiently small, $h<<1/L,$ by (\ref{w20}), and since $y_{0}\in U.$
Then, (\ref{w30-1}) implies (\ref{w19})-(\ref{w19-0}).\hfill $\square $

\subsection{Convergence of the discretization scheme}

Let us define $y^{h}:[0,T]\rightarrow L^{2}(\Omega )$ by
\begin{eqnarray}
y^{h}(t) &=&y_{i}^{h},\mbox{ }t\in \lbrack (i-1)h,ih),\mbox{ }i=1,...n,
\label{w31} \\
y^{h}(0) &=&y_{0}^{h},  \notag
\end{eqnarray}%
and extend $y^{h}$ by continuity to the right of $T$ as
\begin{equation*}
y^{h}(t)=y_{n}^{h},\mbox{ }t\in \lbrack T,T+\delta ],\mbox{ with }\delta
\mbox{ arbitrary, }\delta >h.
\end{equation*}

The step function $y^{h}$ defined by (\ref{w31}) is called an $h$%
-approximating solution to (\ref{w1})-(\ref{w3}) (see \cite{vb-springer-2010}%
, p. 129). Also, we set, for all $i=1,...n,$
\begin{eqnarray*}
f^{h}(t) &=&f_{i}^{h},\mbox{ }t\in \lbrack (i-1)h,ih),\mbox{ } \\
g^{h}(t) &=&g_{i}^{h},\mbox{ }t\in \lbrack (i-1)h,ih),\mbox{ } \\
\beta (t,x,\nabla y^{h}(t)) &=&\beta (t_{i},x,\nabla y_{i}^{h}),\mbox{ }t\in
\lbrack (i-1)h,ih),\mbox{ } \\
\eta ^{h}(t) &=&\eta _{i}^{h},\mbox{ }t\in \lbrack (i-1)h,ih),\mbox{ } \\
j(t,x,\nabla y^{h}(t)) &=&j(t_{i},x,\nabla y_{i}^{h}),\mbox{ }t\in \lbrack
(i-1)h,ih).
\end{eqnarray*}%
We see that $\eta ^{h}(t,x)\in \beta (t,x,\nabla y^{h}(t))$ a.e. on $Q.$

Then, we deduce from (\ref{w17})-(\ref{w20}) the estimates
\begin{equation}
\left\Vert y^{h}(t)\right\Vert _{L^{2}(\Omega )}+\left\Vert \gamma
(y^{h}(t))\right\Vert _{L^{2}(\Gamma )}\leq C,\mbox{ \ for }t\in \lbrack
0,T],  \label{w32}
\end{equation}%
\begin{equation}
\int_{0}^{T}\left\Vert \nabla y^{h}(t)\right\Vert _{L^{p}(\Omega
)}^{p}dt\leq C,  \label{w34}
\end{equation}%
\begin{equation}
\int_{0}^{T}\left\Vert \frac{y^{h}(t+h)-y^{h}(t)}{h}\right\Vert
_{L^{2}(\Omega )}^{2}dt\leq C,  \label{w35}
\end{equation}%
\begin{equation}
\int_{0}^{T}\left\Vert \frac{\gamma (y^{h}(t+h))-\gamma (y^{h}(t))}{h}%
\right\Vert _{L^{2}(\Gamma )}^{2}dt\leq C,  \label{w36}
\end{equation}%
\begin{equation}
\int_{0}^{T}\int_{\Omega }j(t,x,\nabla y^{h}(t))dxdt\leq C,  \label{w37}
\end{equation}%
with $C$ independent of $h.$ Also, (\ref{w32}) and (\ref{w34}) imply that
\begin{equation}
\int_{0}^{T}\left\Vert y^{h}(t)\right\Vert _{W^{1,\overline{p}}(\Omega )}^{%
\overline{p}}dt\leq C,  \label{w37-0}
\end{equation}%
where $\overline{p}$ given by (\ref{w21-5}).

\medskip

\noindent \textbf{Definition 3.1. }We call a \textit{weak} \textit{solution}
to problem (\ref{w1})-(\ref{w3}) a function $y\in L^{2}(Q),$ with
\begin{equation*}
\nabla y\in L^{p}(0,T;(L^{p}(\Omega ))^{N}),\mbox{ }\gamma (y)\in
L^{2}(\Sigma ),\mbox{ }j(\cdot ,\cdot ,\nabla y)\in L^{1}(Q),
\end{equation*}%
\begin{equation*}
\mbox{such that there exists }\eta \in (L^{p^{\prime }}(Q))^{N},\mbox{ }\eta
(t,x)\in \beta (t,x,\nabla y(t,x))\mbox{ a.e. }(t,x)\in Q,\mbox{ satisfying}
\end{equation*}%
\begin{eqnarray}
&&-\int_{Q}y\phi _{t}dxdt+\int_{Q}\eta \cdot \nabla \phi dxdt-\int_{\Sigma
}y\phi _{t}d\sigma dt  \label{w38} \\
&=&\int_{Q}f\phi dxdt+\int_{\Sigma }g\phi d\sigma dt+\int_{\Omega }y_{0}\phi
(0)dx+\int_{\Gamma }y_{0}\phi (0)d\sigma ,  \notag
\end{eqnarray}%
for all $\phi \in W^{1,2}([0,T];L^{2}(\Omega )),$ $\nabla \phi \in
(L^{p}(Q))^{N},$ $\gamma (\phi )\in W^{1,2}([0,T];L^{2}(\Gamma )),$ $\phi
(T)=0.$

\medskip

Theorem 3.2 below is the main result of this section.

\medskip

\noindent \textbf{Theorem 3.2.}\textit{\ Let us assume} (\ref{w21-3})$.$
\textit{Then, under hypotheses }(H$_{1}),$ (H$_{2}),$ (H$_{5}),$ \textit{%
problem }(\ref{w1})-(\ref{w3}) \textit{has at least one weak solution }$y,$
\textit{which satisfies}%
\begin{eqnarray}
y &\in &W^{1,2}([0,T];L^{2}(\Omega )),\mbox{ }\gamma (y)\in
W^{1,2}([0,T];L^{2}(\Gamma )),\mbox{ }  \label{w38-0} \\
\nabla \cdot \eta &\in &L^{2}(Q),\mbox{ }\eta (t,x)\in \beta (t,x,\nabla
y(t,x))\mbox{ \textit{a.e.} }(t,x)\in Q.  \notag
\end{eqnarray}%
\textit{Moreover, }$y$\textit{\ is a strong solution to }(\ref{w1})-(\ref{w3}%
)\textit{, that is}%
\begin{equation}
y_{t}-\nabla \cdot \eta =f,\mbox{ \textit{a.e.} \textit{in} }Q,  \label{w102}
\end{equation}%
\begin{equation}
\gamma (\eta )\cdot \nu +y_{t}=g,\mbox{ \textit{a.e.} \textit{on} }\Sigma ,
\label{w103}
\end{equation}%
\begin{equation}
y(0)=y_{0},\mbox{ \textit{in} }\Omega .  \label{w104}
\end{equation}%
\textit{The solution }$y$ \textit{is given by }
\begin{equation}
y=\lim_{h\rightarrow 0}y^{h}\mbox{\textit{strongly in} }L^{r}(Q),
\label{w105}
\end{equation}%
\textit{with }$y^{h}$\textit{\ defined by }(\ref{w31}),\textit{\ }$r=2$%
\textit{\ if }$p\geq 2$\textit{\ and }$r=p,$\textit{\ if }$p\in (1,2)$%
\textit{. }

\textit{The solution is unique in the class of functions satisfying }(\ref%
{w38-0})-(\ref{w104})\textit{\ and the map }$(y_{0},\gamma
(y_{0}))\rightarrow (y(t),\gamma (y(t)))$ \textit{is Lipschitz from }$%
L^{2}(\Omega )\times L^{2}(\Gamma )$\textit{\ to} $C([0,T];L^{2}(\Omega
))\times C([0,T];L^{2}(\Gamma )).$

\medskip

\noindent \textbf{Proof. }The proof is done in three steps. First, we note
that if $y_{0}\in U$ we have $j(0,\cdot ,\nabla y_{0})\in L^{1}(\Omega ).$

\textbf{Weak solution.} By (\ref{w32})-(\ref{w37}) it follows that one can
select a subsequence such that as $h\rightarrow 0,$ we have,
\begin{equation*}
y^{h}\rightarrow y\mbox{ weak-star in }L^{\infty }(0,T;L^{2}(\Omega )),\mbox{
}
\end{equation*}%
\begin{equation*}
y^{h}\rightarrow y\mbox{ weak-star in }L^{\infty }(0,T;L^{2}(\Gamma )),\mbox{
}
\end{equation*}%
\begin{equation*}
\nabla y^{h}\rightarrow \nabla y\mbox{ weakly in }(L^{p}(Q))^{N},\mbox{ }
\end{equation*}%
\begin{equation*}
\frac{y^{h}(t+h)-y^{h}(t)}{h}\rightarrow l\mbox{ weakly in }L^{2}(Q),\mbox{ }
\end{equation*}%
\begin{equation*}
\frac{\gamma (y^{h}(t+h))-\gamma (y^{h}(t))}{h}\rightarrow l_{1}\mbox{
weakly in }L^{2}(\Sigma ).
\end{equation*}%
To prove the last two assertions we proceed by a direct calculus. For some $%
\delta >0,$ we take $\phi \in M_{\delta },$ where%
\begin{equation*}
M_{\delta }=\{\phi \in C^{\infty }(\overline{Q});\mbox{ }\phi (t,x)=0\mbox{
on }t\in \lbrack T-\delta ,T]\}
\end{equation*}%
and compute (without writing the argument $x$ for the functions $y^{h}$ and $%
\phi )$
\begin{eqnarray*}
&&\int_{0}^{T}\int_{\Omega }\dfrac{y^{h}(t+h)-y^{h}(t)}{h}\phi (t)dxdt \\
&=&\int_{0}^{T-h}\int_{\Omega }\frac{1}{h}y^{h}(t+h)\phi
(t)dxdt-\int_{0}^{T-h}\int_{\Omega }\frac{1}{h}y^{h}(t)\phi (t)dxdt \\
&=&\int_{h}^{T}\int_{\Omega }\frac{1}{h}y^{h}(s)\phi
(s-h)dxds-\int_{0}^{T-h}\int_{\Omega }\frac{1}{h}y^{h}(t)\phi (t)dxdt \\
&=&\int_{h}^{T-h}\int_{\Omega }\frac{1}{h}y^{h}(s)\phi
(s-h)dxds+\int_{T-h}^{T}\int_{\Omega }\frac{1}{h}y^{h}(s)\phi (s-h)dxds \\
&&-\int_{0}^{h}\int_{\Omega }\frac{1}{h}y^{h}(s)\phi
(s)dxds-\int_{h}^{T-h}\int_{\Omega }\frac{1}{h}y^{h}(s)\phi (s)dxds \\
&=&-\int_{h}^{T-h}\int_{\Omega }\frac{\phi (s)-\phi (s-h)}{h}%
y^{h}(s)dxds-\int_{0}^{h}\int_{\Omega }\frac{1}{h}y^{h}(s)\phi (s)dxds \\
&&+\int_{T-h}^{T}\int_{\Omega }\frac{1}{h}y^{h}(s)\phi (s-h)dxds.
\end{eqnarray*}%
In the calculus above $\phi (t,x)=0$ for $t\in \lbrack T-h,T]$ since we can
take $\delta >h.$ Next,
\begin{multline*}
-\int_{0}^{h}\int_{\Omega }\frac{1}{h}y^{h}(s)\phi (s)dxds \\
=-\frac{1}{h}\int_{0}^{h}\int_{\Omega }y^{h}(s)(\phi (s)-\phi (0))dxds-\frac{%
1}{h}\int_{0}^{h}\int_{\Omega }y^{h}(s)\phi (0)dxds \\
\leq \frac{1}{h}\int_{0}^{h}\left\Vert y^{h}(s)\right\Vert _{L^{2}(\Omega
)}\left\Vert \phi (s)-\phi (0)\right\Vert _{L^{2}(\Omega )}ds-\int_{\Omega
}\phi (0)\frac{1}{h}\int_{0}^{h}y^{h}(s)dsdx \\
\leq C\frac{1}{h}\int_{0}^{h}\left\Vert \phi (s)-\phi (0)\right\Vert
_{L^{2}(\Omega )}ds-\int_{\Omega }\phi (0)y_{0}dx+\epsilon (h)\rightarrow
-\int_{\Omega }\phi (0)y_{0}dx,
\end{multline*}%
where $\epsilon (h)\rightarrow 0$ as $h\rightarrow 0$.

Proceeding in the same way for the last term we get
\begin{multline*}
\int_{T-h}^{T}\int_{\Omega }\frac{1}{h}y^{h}(s)\phi
(s-h)dxds=\int_{T-h}^{T}\int_{\Omega }\frac{1}{h}y^{h}(s)(\phi (s-h)-\phi
(T-h))dxds \\
+\int_{T-h}^{T}\int_{\Omega }\frac{1}{h}y^{h}(s)\phi (T-h)dxds\leq C\frac{1}{%
h}\int_{T-h}^{T}\left\Vert \phi (s-h)-\phi (T-h)\right\Vert _{L^{2}(\Omega
)}ds,
\end{multline*}%
as $\phi (T-h)=0.$ Again by the continuity of $\phi $ we obtain that
\begin{equation*}
\int_{T-h}^{T}\int_{\Omega }\frac{1}{h}y^{h}(s)\phi (s-h)dxds\rightarrow 0,%
\mbox{ as }h\rightarrow 0.
\end{equation*}

In conclusion, all these yield
\begin{equation*}
\lim_{h\rightarrow 0}\int_{0}^{T}\int_{\Omega }\dfrac{y^{h}(t+h)-y^{h}(t)}{h}%
\phi (t)dxdt=-\int_{0}^{T}\int_{\Omega }y(t)\frac{d\phi }{dt}%
(t)dxdt-\int_{\Omega }\phi (0)y_{0}dx,
\end{equation*}%
for any $\phi \in M_{\delta }.$ Therefore, we get, in the sense of
distributions, that
\begin{equation}
\lim_{h\rightarrow 0}\int_{0}^{T}\int_{\Omega }\dfrac{y^{h}(t+h)-y^{h}(t)}{h}%
\phi (t)dxdt=\frac{dy}{dt}(\phi ),\mbox{ for any }\phi \in C_{0}^{\infty
}(Q),  \label{w100}
\end{equation}%
and so, $l=\frac{dy}{dt}$ in $\mathcal{D}^{\prime }(Q)$ (the space of
Schwartz distributions on $Q)$. Moreover, by (\ref{w35}) we still have
\begin{equation*}
\left\vert \frac{dy}{dt}(\phi )\right\vert \leq C\left\Vert \phi \right\Vert
_{L^{2}(Q)},\mbox{ for }\phi \in C_{0}^{\infty }(Q)\cap M_{\delta },
\end{equation*}%
with $C$ independent of $\delta ,$ and so $\frac{dy}{dt}\in L^{2}(Q)$ and $l=%
\frac{dy}{dt}$ a.e. on $Q.$ Therefore, $y\in W^{1,2}([0,T-\delta
];L^{2}(\Omega ))$ and since $\delta $ is arbitrary we infer that $y\in
W^{1,2}([0,T];L^{2}(\Omega )).$

Proceeding in the same way for the time differences on $\Gamma $ we get that
\begin{eqnarray}
&&\lim_{h\rightarrow 0}\int_{0}^{T}\int_{\Gamma }\dfrac{\gamma
(y^{h}(t+h))-\gamma (y^{h}(t))}{h}\phi (t)d\sigma dt  \label{w101} \\
&=&-\int_{0}^{T}\int_{\Gamma }\gamma (y(t))\frac{d\phi }{dt}(t)d\sigma
dt-\int_{\Gamma }\phi (0)y_{0}d\sigma ,\mbox{ for }\phi \in M_{\delta }.
\notag
\end{eqnarray}%
Therefore, we obtain that
\begin{equation*}
l_{1}=\frac{d\gamma (y)}{dt}\mbox{ \ a.e. on }\Sigma ,\mbox{ \ }\frac{%
d\gamma (y)}{dt}\in L^{2}(0,T-\delta ;L^{2}(\Gamma ))
\end{equation*}%
and so, finally $\gamma (y)\in W^{1,2}([0,T];L^{2}(\Gamma )).$

We deduce that $\xi =\nabla y$ a.e. on $Q,$ by the same argument used in
Proposition 2.2 and by passing to the limit in (\ref{w37}) and using the
weak lower semicontinuity of the convex integrand we get $j(t,\cdot ,\nabla
y)\in L^{1}(Q).$

The next step is to prove (\ref{w105}). A simple way to show it is to use a
compactness argument in the space of vectorial functions with bounded
variations on $[0,T].$

We have that $y^{h}\in BV([0,T];L^{2}(\Omega )),$ the space of functions
with bounded variation from $[0,T]$ to $L^{2}(\Omega ),$ i.e.,
\begin{equation}
V_{0}^{T}(y^{h})=\sup_{P\in \mathcal{P}}\sum_{i=1}^{n_{p}}\left\Vert
y^{h}(s_{i})-y^{h}(s_{i-1})\right\Vert _{L^{2}(\Omega )}\leq C,  \label{w40}
\end{equation}%
where $C$ is a constant and $\mathcal{P}=\{P=(s_{0},...,s_{n_{p}});P$ is a
partition of $[0,T]\}$ is the set of all partitions of $[0,T].$ Indeed, if
we consider an equidistant partition (e.g., with $s_{i}=t_{i}$) we have by (%
\ref{w19}) that
\begin{eqnarray}
&&\left( \sum_{i=0}^{n-1}\left\Vert y^{h}(t_{i+1})-y^{h}(t_{i})\right\Vert
_{L^{2}(\Omega )}\right) ^{2}  \label{w41} \\
&\leq &n\sum_{i=0}^{n-1}\left\Vert y_{i+1}^{h}-y_{i}^{h}\right\Vert
_{L^{2}(\Omega )}^{2}=nh\cdot h\sum_{i=1}^{n}\left\Vert \frac{%
y_{i+1}^{h}-y_{i}^{h}}{h}\right\Vert _{L^{2}(\Omega )}^{2}\leq TC.  \notag
\end{eqnarray}

Now, we discuss separately the cases $p\geq 2$ and $p\in (1,2).$

Let $p\geq 2.$ By (\ref{w41}) we also have $y^{h}\in BV([0,T];(H^{1}(\Omega
))^{\prime }).$

On the basis of this relation, (\ref{w32}), and since $L^{2}(\Omega )$ is
compact in $(H^{1}(\Omega ))^{\prime }$ we can apply the strong version of
Helly theorem for the infinite dimensional case (see \cite{vbp-2012}, Remark
1.127, p. 48). We deduce that
\begin{equation}
y^{h}(t)\rightarrow y(t)\mbox{ strongly in }(H^{1}(\Omega ))^{\prime }\mbox{
uniformly in }t\in \lbrack 0,T].  \label{w42}
\end{equation}

Next, applying Lemma\ 5.1 in \cite{lions}, p. 58, we have that for any $%
\varepsilon >0$ there exists a constant $C_{\varepsilon }$ such that
\begin{equation}
\left\Vert w\right\Vert _{L^{2}(\Omega )}\leq \varepsilon \left\Vert
w\right\Vert _{H^{1}(\Omega )}+C_{\varepsilon }\left\Vert w\right\Vert
_{(H^{1}(\Omega ))^{\prime }},\mbox{ }\forall w\in H^{1}(\Omega ).
\label{w43}
\end{equation}%
This lemma\ applied for $w=y^{h}(t)-y(t)$ yields
\begin{equation*}
\frac{1}{2}\left\Vert y^{h}(t)-y(t)\right\Vert _{L^{2}(\Omega )}^{2}\leq
\varepsilon \left\Vert y^{h}(t)-y(t)\right\Vert _{H^{1}(\Omega
)}^{2}+C_{\varepsilon }\left\Vert y^{h}(t)-y(t)\right\Vert _{(H^{1}(\Omega
))^{\prime }}^{2}.
\end{equation*}%
Integrating with respect to $t$ on $(0,T)$ we obtain that
\begin{eqnarray*}
&&\frac{1}{2}\int_{0}^{T}\left\Vert y^{h}(t)-y(t)\right\Vert _{L^{2}(\Omega
)}^{2}dt \\
&\leq &\varepsilon \int_{0}^{T}\left\Vert y^{h}(t)-y(t)\right\Vert
_{H^{1}(\Omega )}^{2}dt+C_{\varepsilon }\int_{0}^{T}\left\Vert
y^{h}(t)-y(t)\right\Vert _{(H^{1}(\Omega ))^{\prime }}^{2}dt.
\end{eqnarray*}%
Then, the last term on the right-hand side tends to 0 as $h\rightarrow 0,$
by (\ref{w42}), and the coefficient of $\varepsilon $ is bounded, by (\ref%
{w37-0}). Hence
\begin{equation*}
\limsup_{h\rightarrow 0}\int_{0}^{T}\left\Vert y^{h}(t)-y(t)\right\Vert
_{L^{2}(\Omega )}^{2}dt\leq C\varepsilon \mbox{ for any }\varepsilon >0,
\end{equation*}%
and since $\varepsilon $ is arbitrary we get (\ref{w105}), with $r=2.$

Let $p\in (1,2).$ We shall prove that
\begin{equation}
y^{h}\rightarrow y\mbox{ strongly in }L^{p}(Q),\mbox{ as }h\rightarrow 0.
\label{w42-0}
\end{equation}%
From (\ref{w41}) it follows that $y^{h}\in BV([0,T];L^{p}(\Omega )).$ We
assert that there exists a Banach space $X$ such that $L^{p}(\Omega )\subset
X,$ with compact injection. For example, $X=(W^{1,r}(\Omega ))^{\prime }$
(i.e., if $W^{1,r}(\Omega )\subset L^{p^{\prime }}(\Omega ),$ with $%
p^{\prime }<\frac{rN}{N-r}$ if $N>r,$ and $p^{\prime }=r$ if $N\leq r).$

\noindent Therefore, once again by the strong Helly theorem and since $%
\left\Vert y^{h}(t)\right\Vert _{L^{p}(\Omega )}\leq C$ we get
\begin{equation}
y^{h}(t)\rightarrow y(t)\mbox{ strongly in }X,\mbox{ uniformly with }t\in
\lbrack 0,T].  \label{w42-1}
\end{equation}

Then, we apply the argument before for the triplet $W^{1,p}(\Omega )\subset
L^{p}(\Omega )\subset X$ and use (\ref{w42-1}) and (\ref{w37-0}). We get
\begin{eqnarray*}
&&\frac{1}{2}\int_{0}^{T}\left\Vert y^{h}(t)-y(t)\right\Vert _{L^{p}(\Omega
)}^{p}dt \\
&\leq &\varepsilon \int_{0}^{T}\left\Vert y^{h}(t)-y(t)\right\Vert
_{W^{1,p}(\Omega )}^{p}dt+C_{\varepsilon }\int_{0}^{T}\left\Vert
y^{h}(t)-y(t)\right\Vert _{X}^{p}dt\rightarrow 0,\mbox{ }
\end{eqnarray*}%
as $\varepsilon \rightarrow 0,$ whence (\ref{w42-0}) follows.

Let us fix $t\in \lbrack 0,T].$ By (\ref{w32}), we have that on a
subsequence,
\begin{equation*}
y^{h}(t)\rightarrow \vartheta (t)\mbox{ weakly in }L^{2}(\Omega ),\mbox{ as }%
h\rightarrow 0,\mbox{ for each fixed }t\in \lbrack 0,T],
\end{equation*}%
and since we have either (\ref{w42}) or (\ref{w42-1}) we get by the limit
uniqueness that $\vartheta (t)=y(t),$ a.e. on $\Omega .$ In particular, it
follows that
\begin{equation*}
y^{h}(T)\rightarrow y(T)\mbox{ weakly in }L^{2}(\Omega ),\mbox{ as }%
h\rightarrow 0.
\end{equation*}

We mention, that as $h\rightarrow 0,$
\begin{equation}
f^{h}\rightarrow f\mbox{ strongly in }L^{2}(0,T;L^{2}(\Omega )),  \label{w45}
\end{equation}%
\begin{equation}
g^{h}\rightarrow g\mbox{ strongly in }L^{2}(0,T;L^{2}(\Gamma )).
\label{w45-1}
\end{equation}

Now, relation (\ref{w21-2}) yields
\begin{equation*}
\left\vert \eta ^{h}(t,x)\right\vert \leq C_{3}\left\vert \nabla
y^{h}(t,x)\right\vert ^{p-1}+C_{3}^{0},
\end{equation*}%
for $\eta ^{h}(t,x)\in \beta (t,x,\nabla y^{h}(t,x))$ a.e. on $Q,$ and so we
conclude that $\{\eta ^{h}\}_{h}$ is bounded in $(L^{p^{\prime }}(Q))^{N}.$
Therefore, on a subsequence, we have
\begin{equation*}
\eta ^{h}\rightarrow \eta \mbox{ weakly in }(L^{p^{\prime }}(Q))^{N},\mbox{
as }h\rightarrow 0,\mbox{ }
\end{equation*}%
and it remains to prove that $\eta (t,x)\in \beta (t,x,\nabla y(t,x)),$ a.e.
$(t,x)\in Q.$

Summing up (\ref{w24-0}) from $i=0$\ to $n-1$ we get
\begin{multline*}
\frac{1}{2}\left\Vert y_{n}^{h}\right\Vert _{L^{2}(\Omega
)}^{2}+h\sum_{i=0}^{n-1}\int_{\Omega }\eta _{i+1}^{h}\cdot \nabla
y_{i+1}^{h}dx+\frac{1}{2}\left\Vert y_{n}^{h}\right\Vert _{L^{2}(\Gamma
)}^{2} \\
\leq h\sum_{i=0}^{n-1}\int_{\Omega
}f_{i+1}^{h}y_{i+1}^{h}dx+h\sum_{i=0}^{n-1}\int_{\Gamma
}g_{i+1}^{h}y_{i+1}^{h}d\sigma +\frac{1}{2}\left\Vert y_{0}\right\Vert
_{L^{2}(\Omega )}^{2}+\frac{1}{2}\left\Vert y_{0}\right\Vert _{L^{2}(\Gamma
)}^{2},
\end{multline*}%
whence, replacing the definition of $y^{h}(t)$ we obtain
\begin{eqnarray*}
&&\frac{1}{2}\left\Vert y^{h}(T)\right\Vert _{L^{2}(\Omega
)}^{2}+h\int_{0}^{T}\int_{\Omega }\eta ^{h}\cdot \nabla y^{h}dxdt+\frac{1}{2}%
\left\Vert y^{h}(T)\right\Vert _{L^{2}(\Gamma )}^{2} \\
&\leq &h\int_{0}^{T}\int_{\Omega }f^{h}y^{h}dxdt+h\int_{0}^{T}\int_{\Gamma
}g^{h}y^{h}d\sigma dt+\frac{1}{2}\left\Vert y_{0}\right\Vert _{L^{2}(\Omega
)}^{2}+\frac{1}{2}\left\Vert y_{0}\right\Vert _{L^{2}(\Gamma )}^{2}.
\end{eqnarray*}%
Passing to the limit as $h\rightarrow 0,$ this yields
\begin{eqnarray}
&&\limsup\limits_{h\rightarrow 0}\int_{Q}\eta ^{h}\cdot \nabla y^{h}dxdt\leq
-\frac{1}{2}\left\Vert y(T)\right\Vert _{L^{2}(\Omega )}^{2}  \label{w47-2}
\\
&&+\int_{Q}fydxdt+\int_{\Sigma }gyd\sigma dt+\frac{1}{2}\left\Vert
y_{0}\right\Vert _{L^{2}(\Omega )}^{2}+\frac{1}{2}\left\Vert
y_{0}\right\Vert _{L^{2}(\Gamma )}^{2}-\frac{1}{2}\left\Vert y(T)\right\Vert
_{L^{2}(\Gamma )}^{2}.  \notag
\end{eqnarray}

Now we write (\ref{w12})-(\ref{w13}) in the following form, after replacing
the functions $y_{i}^{h}$ by $y^{h},$ and integrating with respect to $t\in
(0,T)$
\begin{multline}
\int_{Q}\frac{y^{h}(t+h,x)-y^{h}(t,x)}{h}\phi (t,x)dxdt+\int_{Q}\eta
^{h}(t+h,x)\cdot \nabla \phi (t,x)dxdt  \label{w47-0} \\
+\int_{\Sigma }\frac{\gamma (y^{h}(t+h,x))-\gamma (y^{h}(t,x))}{h}\phi
(t,x)d\sigma dt  \notag \\
=\int_{Q}f^{h}(t+h,x)\phi (t,x)dxdt+\int_{\Sigma }g^{h}(t+h,x)\phi
(t,x)d\sigma dt  \notag
\end{multline}%
for any $\phi \in C^{\infty }(\overline{Q}),$ with $\phi (T,x)=0$. We pass
to the limit as $h\rightarrow 0,$ using the convergences previously deduced
and get
\begin{eqnarray}
&&-\int_{Q}y\phi _{t}dxdt+\int_{Q}\eta \cdot \nabla \phi dxdt-\int_{\Sigma
}y\phi _{t}d\sigma dt  \label{w47-1} \\
&=&\int_{Q}f\phi dxdt+\int_{\Sigma }g\phi d\sigma dt+\int_{\Omega }y_{0}\phi
(0)dx+\int_{\Gamma }y_{0}\phi (0)d\sigma ,  \notag
\end{eqnarray}%
where $\eta =\lim\limits_{h\rightarrow 0}\eta ^{h}$ (weakly in $%
(L^{p^{\prime }}(\Omega ))^{N}$).

Taking into account that $y\in W^{1,2}([0,T];L^{2}(\Omega ))$ we have
\begin{equation*}
\int_{Q}y\phi _{t}dxdt=-\int_{Q}y_{t}\phi dxdt-\int_{\Omega }y_{0}\phi (0)dx
\end{equation*}%
and so we obtain%
\begin{equation}
\int_{Q}y_{t}\phi dxdt+\int_{Q}\eta \cdot \nabla \phi dxdt+\int_{\Sigma
}y_{t}\phi d\sigma dt=\int_{Q}f\phi dxdt+\int_{\Sigma }g\phi d\sigma dt.
\label{w47-3}
\end{equation}

By density, this extends to all $\phi \in W^{1,2}([0,T];L^{2}(\Omega ))$
with $\nabla \phi \in (L^{p}(Q))^{N},$ and in particular, for $\phi =y.$
Finally, we have got%
\begin{eqnarray}
\int_{Q}\eta \cdot \nabla ydxdt &=&-\frac{1}{2}\left\Vert y(T)\right\Vert
_{L^{2}(\Omega )}^{2}+\frac{1}{2}\left\Vert y_{0}\right\Vert _{L^{2}(\Omega
)}^{2}  \label{w47-4} \\[1pt]
&&+\int_{Q}fydxdt+\int_{\Sigma }g\phi d\sigma dt-\frac{1}{2}\left\Vert
y(T)\right\Vert _{L^{2}(\Gamma )}^{2}+\frac{1}{2}\left\Vert y_{0}\right\Vert
_{L^{2}(\Gamma )}^{2}.  \notag
\end{eqnarray}%
Comparing with (\ref{w47-2}) we deduce that
\begin{equation*}
\limsup\limits_{h\rightarrow 0}\int_{0}^{T}\int_{\Omega }\eta ^{h}\cdot
\nabla y^{h}dxdt\leq \int_{Q}\eta \cdot \nabla ydxdt
\end{equation*}%
and since the operator $z\rightarrow \beta (t,x,z)$ is maximal monotone in
the dual pair $(L^{p}(Q))^{N}$ - $L^{p^{\prime }}(Q))^{N},$ we get $\eta
(t,x)\in \beta (t,x,\nabla y(t,x)),$ a.e. on $Q$ (see e.g., \cite%
{vb-springer-2010}, p. 41). Hence $y$ is a weak solution to (\ref{w1})-(\ref%
{w3}).

\medskip

\textbf{Strong solution.} By (\ref{w47-3}) we see that
\begin{equation*}
\int_{Q}y_{t}\phi dxdt+\int_{Q}\eta \cdot \nabla \phi dxdt=\int_{Q}f\phi
dxdt,\mbox{ }\forall \phi \in C_{0}^{\infty }(Q).
\end{equation*}%
Since, we have
\begin{equation*}
(-\nabla \cdot \eta )(\phi )=\int_{Q}\eta \cdot \nabla \phi dxdt,\mbox{ }%
\forall \phi \in C_{0}^{\infty }(Q),
\end{equation*}%
and recalling that $y_{t},$ $f\in L^{2}(Q),$ $\eta \in (L^{p^{\prime
}}(Q))^{N},$ we get (\ref{w102}) in $\mathcal{D}^{\prime }(Q),$ as claimed.

We note that since $\nabla \cdot \eta \in L^{2}(Q)$ it follows that
\begin{equation*}
\gamma (\eta )\cdot \nu \in L^{2}(0,T;H^{-1/2}(\Gamma ))
\end{equation*}%
is well-defined (see e.g., \cite{andreu}, or Theorem 1.2. in \cite{temam})
and the following formula holds
\begin{equation}
\int_{\Omega }\phi (t)\nabla \cdot \eta (t)dx=-\int_{\Omega }\eta (t)\cdot
\nabla \phi (t)dx+\int_{\Gamma }\phi (t)(\gamma (\eta (t))\cdot \nu )d\sigma
,\mbox{ a.e. }t\in (0,T).  \label{w47-4-1}
\end{equation}

Next, we multiply (\ref{w102}) by $\phi \in W^{1,2}([0,T];L^{2}(\Omega ))$
with $\nabla \phi \in (L^{p}(Q))^{N},$ $\gamma (\phi )\in
W^{1,2}([0,T];L^{2}(\Gamma ))$, $\phi (T)=0$ and by (\ref{w47-4-1}) we get
that
\begin{equation}
-\int_{Q}y\phi _{t}dxdt-\int_{\Omega }y_{0}\phi (0)dx-\int_{\Sigma }\phi
\eta \cdot \nu d\sigma dt+\int_{Q}\eta \cdot \nabla \phi dxdt=\int_{Q}f\phi
dxdt.  \label{w47-5}
\end{equation}%
After replacing (\ref{w47-5}) in (\ref{w38}), we obtain that
\begin{equation*}
\int_{\Sigma }\gamma (y)\phi _{t}d\sigma dt+\int_{\Gamma }y_{0}\phi
(0)d\sigma -\int_{\Sigma }\phi \eta \cdot \nu d\sigma dt+\int_{\Sigma }g\phi
d\sigma dt=0,
\end{equation*}%
whence we obtain (\ref{w103}) in the sense of distributions and also a.e. on
$\Sigma ,$ since, as seen earlier, $\frac{d}{dt}\gamma (y)\in L^{2}(\Sigma
). $

\medskip

\textbf{Continuous dependence on data. }Let us consider two solutions $y$
and $\overline{y}$ to (\ref{w1})-(\ref{w3}), corresponding to the data $%
(y_{0},f,g)$ and $(\overline{y_{0}},\overline{f},\overline{g}),$
respectively, in the class of functions satisfying (\ref{w38-0})-(\ref{w104}%
). We make the difference of the two equations (\ref{w102}), corresponding
to these data, multiply the difference by $(y-\overline{y})(t)$ and
integrate on $\Omega .$ We get, a.e. $t\in \lbrack 0,T]$
\begin{multline*}
\frac{1}{2}\frac{d}{dt}\left\Vert y(t)-\overline{y}(t)\right\Vert
_{L^{2}(\Omega )}^{2}+\int_{\Omega }(\eta (t)-\overline{\eta }(t))\cdot
(\nabla y(t)-\nabla \overline{y}(t))dx \\
+\frac{1}{2}\frac{d}{dt}\left\Vert \gamma (y(t))-\gamma (\overline{y}%
(t))\right\Vert _{L^{2}(\Gamma )}^{2}=\int_{\Omega }(f(t)-\overline{f}%
(t))(y(t)-\overline{y}(t))dt,
\end{multline*}%
where $\eta (t,x)\in \beta (t,x,\nabla y(t,x)),$ $\overline{\eta }(t,x)\in
\beta (t,x,\nabla \overline{y}(t,x))$ a.e. $(t,x)\in Q.$ Using the
monotonicity of $\beta $ and integrating with respect to $t$ we get that
\begin{eqnarray*}
&&\left\Vert y(t)-\overline{y}(t)\right\Vert _{L^{2}(\Omega
)}^{2}+\left\Vert \gamma (y(t))-\gamma (\overline{y}(t))\right\Vert
_{L^{2}(\Gamma )}^{2} \\
&\leq &C\left(
\begin{array}{c}
\left\Vert y_{0}-\overline{y_{0}}\right\Vert _{L^{2}(\Omega
)}^{2}+\left\Vert \gamma (y_{0})-\gamma (\overline{y_{0}})\right\Vert
_{L^{2}(\Gamma )}^{2} \\
+\int_{0}^{T}\left\Vert f(t)-\overline{f}(t)\right\Vert _{L^{2}(\Omega
)}^{2}dt+\int_{0}^{T}\left\Vert g(t)-\overline{g}(t)\right\Vert
_{L^{2}(\Gamma )}^{2}dt%
\end{array}%
\right) ,\mbox{ }\forall t\in \lbrack 0,T],
\end{eqnarray*}%
as claimed.\hfill $\square $

\section{The weakly coercive case}

\setcounter{equation}{0}

In this section we assume that hypotheses (H$_{1})$, (H$_{3})$-(H$_{5})$
hold. Also, as mentioned earlier, without loss of generality we may assume (%
\ref{w4-0}).

We note that here no polynomial growth or coercivity on $j$ are assumed
whatever.

A standard example in this case is
\begin{equation*}
\beta (t,x,r)=a(t,x)\log (\left\vert r\right\vert _{N}+1)\mbox{sgn }r+a(t,x)%
\frac{r}{\left\vert r\right\vert _{N}+1},
\end{equation*}%
where $a\in C^{1}(\overline{Q}),$ $a>0,$ sgn $r=r/\left\vert r\right\vert
_{N}$ for $r\neq 0,$ sgn $0=\{r;$ $\left\vert r\right\vert _{N}\leq 1\}.$

\medskip

On the other hand, monotone functions $r\rightarrow \beta (t,x,r)$ with
exponential growth and symmetric at $\pm \infty ,$ in the sense of (\ref%
{w296}), are accepted by the current hypotheses.

First, we note down for later use the following simple lemma.

\medskip

\noindent \textbf{Lemma 4.1.}\textit{\ Let }%
\begin{equation*}
u\in (L^{1}(\Omega ))^{N},\mathit{\ }w\in (L^{1}(\Omega ))^{N},\mbox{ }%
j(\cdot ,\cdot ,u)\in L^{1}(\Omega ),\mbox{ }j^{\ast }(\cdot ,\cdot ,w)\in
L^{1}(\Omega ).
\end{equation*}%
\textit{\ Then, under assumption }(\ref{w296}) \textit{we have}
\begin{equation}
u\cdot w\in L^{1}(\Omega ).  \label{w326-0}
\end{equation}

\medskip

\noindent \textbf{Proof. }First, we recall the relations (see e.g., \cite%
{vb-springer-2010}, p. 8)
\begin{equation}
j(t,x,r)+j^{\ast }(t,x,\omega )\geq \omega \cdot r,\mbox{ }\forall r,\omega
\in \mathbb{R}^{N}\mbox{, }\forall (t,x)\in \overline{Q}  \label{w291}
\end{equation}%
\begin{equation}
j(t,x,r)+j^{\ast }(t,x,\omega )=\omega \cdot r\mbox{ iff }\omega \in
\partial j(t,x,r),\mbox{ }\forall (t,x)\in \overline{Q}.  \label{w292}
\end{equation}%
By (\ref{w291}),
\begin{equation*}
j(t,x,u(x))+j^{\ast }(t,x,w(x))\geq u(x)\cdot w(x),\mbox{ }\forall (t,x)\in
\overline{Q},
\end{equation*}%
and this yields
\begin{equation*}
\int_{\Omega }u(x)\cdot w(x)dx<\infty .
\end{equation*}%
We write (\ref{w291}) for $(-u^{\ast })$
\begin{equation*}
j(t,x,-u^{\ast }(x))+j^{\ast }(t,x,w(x))\geq -u^{\ast }(x)\cdot w(x),\mbox{ }%
\forall (t,x)\in \overline{Q},
\end{equation*}%
and use (\ref{w296}), obtaining
\begin{equation*}
\int_{\Omega }\left( -u\cdot w\right) dx\leq \gamma _{1}\int_{\Omega
}j(t,x,u)dx+\gamma _{2}\mbox{meas}(\Omega )+\int_{\Omega }j^{\ast
}(t,x,w)dx<\infty .
\end{equation*}%
Therefore we get (\ref{w326-0}), as claimed.\hfill $\square $

\medskip

Let us define the space%
\begin{equation}
U_{1}=\{z\in L^{2}(\Omega );\mbox{ }z\in W^{1,1}(\Omega ),\mbox{ }\gamma
(z)\in L^{2}(\Gamma )\}.  \label{w304}
\end{equation}

For $t$ fixed in $[0,T]$, $h>0,$ and $w_{1}\in L^{2}(\Omega ),$ $w_{2}\in
L^{2}(\Gamma ),$ let us consider the problem%
\begin{eqnarray}
u-h\nabla \cdot \beta (t,x,\nabla u) &\ni &w_{1}\mbox{ in }\Omega ,
\label{w200} \\
u+h\beta (t,x,\nabla u)\cdot \nu &\ni &w_{2}\mbox{ on }\Gamma .  \notag
\end{eqnarray}

As in the previous case, we call a \textit{weak solution} to problem (\ref%
{w200}) a function $u\in U_{1},$ such that $j(t,\cdot ,\nabla u)\in
L^{1}(\Omega ),$ and there exists
\begin{equation}
\eta \in (L^{1}(\Omega ))^{N},\mbox{ }\eta (x)\in \beta (t,x,\nabla u(x))%
\mbox{ a.e. }x\in \Omega ,\mbox{ }j^{\ast }(t,\cdot ,\eta )\in L^{1}(\Omega
),  \label{w200-0}
\end{equation}%
satisfying%
\begin{equation}
\int_{\Omega }(u\psi +h\eta \cdot \nabla \psi )dx+\int_{\Gamma }u\psi
d\sigma =b(\psi ),\mbox{ }\forall \psi \in C^{1}(\overline{\Omega }),
\label{w301}
\end{equation}%
with $b$ given by (\ref{w20-00}) for all $u\in U_{1}.$

Problem (\ref{w200}) has a unique solution, namely given by the unique
minimizer of the functional $\varphi :L^{2}(\Omega )\rightarrow \mathbb{R},$%
\begin{equation}
\varphi (u)=\left\{
\begin{array}{l}
\frac{1}{2}\int_{\Omega }u^{2}dx+h\int_{\Omega }j(t,x,\nabla u)dx+\frac{1}{2}%
\int_{\Gamma }u^{2}dx-b(u),\mbox{ if }u\in U_{1}, \\
+\infty ,\mbox{ \ \ \ \ \ \ \ \ \ \ \ \ \ \ \ \ \ \ \ \ \ \ \ \ \ \ \ \ \ \
\ \ \ \ \ \ \ \ \ \ \ \ \ \ \ \ \ \ \ \ \ \ \ \ \ \ \ \ otherwise.}%
\end{array}%
\right.  \label{w303}
\end{equation}

Actually, we have the equivalence between (\ref{w200}) and the minimization
problem
\begin{equation}
\mbox{Min }\left\{ \varphi (u);\mbox{ }u\in U_{1}\right\} .  \label{w302}
\end{equation}

\medskip

\noindent \textbf{Proposition 4.2. }\textit{Problem}\textbf{\textit{\ }}(\ref%
{w200})\textbf{\ }\textit{has a unique solution which is the minimizer of }$%
\varphi .$

\medskip

\noindent \textbf{Proof. }Let $\lambda >0$ and consider the approximating
regularized problem
\begin{eqnarray}
u-h\nabla \cdot (\beta _{\lambda }(t,x,\nabla u)+\lambda \nabla u) &=&w_{1}%
\mbox{ in }\Omega ,  \label{w305} \\
u+h(\beta _{\lambda }(t,x,\nabla u)+\lambda \nabla u)\cdot \nu &=&w_{2}\mbox{
on }\Gamma ,  \notag
\end{eqnarray}%
where $\beta _{\lambda }$ is the Yosida approximation of $\beta ,$
\begin{equation}
\beta _{\lambda }(t,x,r)=\frac{1}{\lambda }(1-(1+\lambda \beta (t,x,\cdot
))^{-1})r,\mbox{ }\forall r\in \mathbb{R}^{N}.  \label{w306}
\end{equation}%
Its potential (i.e., the Moreau regularization of $j)$ is given by
\begin{eqnarray}
j_{\lambda }(t,x,r) &=&\inf_{s\in \mathbb{R}^{N}}\left\{ \frac{\left\vert
r-s\right\vert _{N}^{2}}{2\lambda }+j(t,x,s)\right\}  \label{w307} \\
&=&\frac{1}{2\lambda }\left\vert (1+\lambda \beta (t,x,\cdot
))^{-1}r-r\right\vert _{N}^{2}+j(t,x,(1+\lambda \beta (t,x,\cdot ))^{-1}r),
\notag
\end{eqnarray}%
and the function $j_{\lambda }$ has the following properties%
\begin{eqnarray}
j_{\lambda }(t,x,r) &\leq &j(t,x,r)\mbox{ for all }r\in \mathbb{R}^{N},\mbox{
}(t,x)\in \overline{Q},\mbox{ }\lambda >0,  \label{w307-0} \\
\lim_{\lambda \rightarrow 0}j_{\lambda }(t,x,r) &=&j(t,x,r),\mbox{ for all }%
r\in \mathbb{R}^{N},\mbox{ }(t,x)\in \overline{Q}.  \notag
\end{eqnarray}

As in Proposition 2.2 we deduce that the solution to (\ref{w305}) is
provided by the unique minimizer of the problem
\begin{equation}
\mbox{Min }\left\{ \varphi _{\lambda }(u);\mbox{ }u\in L^{2}(\Omega
)\right\} ,  \label{w308}
\end{equation}%
where $\varphi _{\lambda }:L^{2}(\Omega )\rightarrow \mathbb{R},$
\begin{equation}
\varphi _{\lambda }(u)=\left\{
\begin{array}{l}
\frac{1}{2}\int_{\Omega }u^{2}dx+h\int_{\Omega }j_{\lambda }(t,x,\nabla u)dx+%
\frac{1}{2}\int_{\Gamma }u^{2}d\sigma +\lambda \int_{\Omega }\left\vert
\nabla u\right\vert _{N}^{2}dx-b(u), \\
\mbox{ \ \ \ \ \ \ \ \ \ \ \ \ \ \ \ \ \ \ \ \ \ \ \ \ \ \ \ \ \ \ \ \ \ \ \
\ \ \ \ \ \ \ \ \ \ \ \ \ \ \ \ \ \ \ \ if }u\in W^{1,2}(\Omega ), \\
+\infty ,\mbox{ \ \ \ \ \ \ \ \ \ \ \ \ \ \ \ \ \ \ \ \ \ \ \ \ \ \ \ \ \ \
\ \ \ \ \ \ \ \ \ \ \ \ \ \ \ \ \ \ otherwise.}%
\end{array}%
\right.  \label{w309}
\end{equation}

Namely, we have, following Proposition 2.2, that in this case the weak
solution $u_{\lambda }\in W^{1,2}(\Omega )$ and it satisfies
\begin{eqnarray}
&&\int_{\Omega }u_{\lambda }\psi dx+h\int_{\Omega }\beta _{\lambda
}(t,x,\nabla u_{\lambda })\cdot \nabla \psi dx+\lambda \int_{\Omega }\nabla
u_{\lambda }\cdot \nabla \psi dx+\int_{\Gamma }u_{\lambda }\psi d\sigma
\notag \\
&=&\int_{\Omega }w_{1}\psi dx+\int_{\Gamma }w_{2}\psi d\sigma ,\mbox{ }%
\forall \psi \in W^{1,2}(\Omega ).  \label{w310}
\end{eqnarray}%
In particular for $\psi =u_{\lambda },$ this yields%
\begin{eqnarray}
&&\frac{1}{2}\int_{\Omega }u_{\lambda }^{2}dx+\int_{\Omega }\beta _{\lambda
}(t,x,\nabla u_{\lambda })\cdot \nabla u_{\lambda }dx+\lambda h\int_{\Omega
}\left\vert \nabla u_{\lambda }\right\vert _{N}^{2}dx+\int_{\Gamma
}u_{\lambda }^{2}d\sigma  \notag \\
&=&\int_{\Omega }w_{1}u_{\lambda }dx+\int_{\Gamma }w_{2}u_{\lambda }d\sigma ,
\label{w312}
\end{eqnarray}%
whence we obtain the estimate
\begin{equation}
\int_{\Omega }u_{\lambda }^{2}dx+2h\int_{\Omega }j_{\lambda }(t,x,\nabla
u_{\lambda })dx+2\lambda h\int_{\Omega }\left\vert \nabla u_{\lambda
}\right\vert _{N}^{2}dx+\int_{\Gamma }u_{\lambda }^{2}d\sigma \leq C.
\label{w313}
\end{equation}%
(By $C$ we denote a positive constant independent of $\lambda .)$

Replacing the definition (\ref{w307}) (second line) of $j_{\lambda }$ we get
\begin{eqnarray}
&&\int_{\Omega }u_{\lambda }^{2}dx+h\int_{\Omega }\frac{1}{\lambda }%
\left\vert (1+\lambda \beta (t,x,\cdot ))^{-1}\nabla u_{\lambda }-\nabla
u_{\lambda }\right\vert _{N}^{2}dx  \label{w314} \\
&&+2h\int_{\Omega }j(t,x,(1+\lambda \beta (t,x,\cdot ))^{-1}\nabla
u_{\lambda })dx+2\lambda \int_{\Omega }\left\vert \nabla u_{\lambda
}\right\vert _{N}^{2}dx+\int_{\Gamma }u_{\lambda }^{2}d\sigma \leq C.  \notag
\end{eqnarray}%
Consequently, each term on the left-hand side in (\ref{w314}) is bounded
independently of $\lambda $ and in particular
\begin{equation}
\int_{\Omega }j(t,x,(1+\lambda \beta (t,x,\cdot ))^{-1}\nabla u_{\lambda
})dx\leq C,\mbox{ }\forall \lambda >0.  \label{w322}
\end{equation}

By (\ref{w322}), (\ref{w293}) and the Dunford-Pettis theorem we can deduce
that the sequence $\left\{ (1+\lambda \beta (t,x,\cdot ))^{-1}\nabla
u_{\lambda }\right\} _{\lambda >0}$ is weakly compact in $L^{1}(\Omega ).$
Indeed, denoting $z_{\lambda }=(1+\lambda \beta (t,x,\cdot ))^{-1}\nabla
u_{\lambda }$, we have to show that the integrals $\int_{S}\left\vert
z_{\lambda }\right\vert _{N}dx,$ with $S\subset \Omega ,$ are
equi-absolutely continuous, meaning that for every $\varepsilon >0$ there
exists $\delta $ such that $\int_{S}\left\vert z_{\lambda }\right\vert
_{N}dxdt<\varepsilon $ whenever meas$(S)<\delta .$ Let $M_{\varepsilon }>%
\frac{2C}{\varepsilon },$ where $C$ is the constant in (\ref{w322}), and let
$R_{M}$ be such that $\frac{j(t,x,z_{\lambda })}{\left\vert z_{\lambda
}\right\vert _{N}}\geq M_{\varepsilon }$ for $\left\vert z_{\lambda
}\right\vert _{N}>R_{M},$ by (\ref{w293}). If $\delta <\frac{\varepsilon }{%
2R_{M}}$ then
\begin{eqnarray*}
&&\int_{S}\left\vert z_{\lambda }\right\vert dxdt\leq \int_{\{x;\left\vert
z_{\lambda }(x)\right\vert _{N}\geq R_{M}\}}\left\vert z_{\lambda
}\right\vert dx+\int_{\{x;\left\vert z_{\lambda }(x)\right\vert
_{N}<R_{M}\}}\left\vert z_{\lambda }\right\vert dx \\
&\leq &M_{\varepsilon }^{-1}\int_{\Omega }j(t,x,z_{\lambda
}(x))dx+R_{M}\delta <\varepsilon .
\end{eqnarray*}%
Then, we select a subsequence (again denoted $_{\lambda })$ such that as $%
\lambda \rightarrow 0,$ we have
\begin{equation}
u_{\lambda }\rightarrow u^{\ast }\mbox{ weakly in }L^{2}(\Omega ),\mbox{ }
\label{w315}
\end{equation}%
\begin{equation}
\gamma (u_{\lambda })\rightarrow \gamma (u^{\ast })\mbox{ weakly in }%
L^{2}(\Gamma ),\mbox{ }  \label{w315-0}
\end{equation}%
\begin{equation}
(1+\lambda \beta (t,x,\cdot ))^{-1}\nabla u_{\lambda }\rightarrow \zeta _{1}%
\mbox{ weakly in }(L^{1}(\Omega ))^{N},\mbox{ }  \label{w316}
\end{equation}%
\begin{equation}
(1+\lambda \beta (t,x,\cdot ))^{-1}\nabla u_{\lambda }-\nabla u_{\lambda
}\rightarrow 0\mbox{ strongly in }(L^{2}(\Omega ))^{N},\mbox{ }  \label{w317}
\end{equation}%
\begin{equation}
\sqrt{\lambda }\nabla u_{\lambda }\rightarrow \zeta _{2}\mbox{ weakly in }%
(L^{2}(\Omega ))^{N}.  \label{w318}
\end{equation}%
By (\ref{w316}), (\ref{w317}) and (\ref{w315}) we get that
\begin{equation}
\nabla u_{\lambda }\rightarrow \zeta _{1}=\nabla u^{\ast }\mbox{ weakly in }%
(L^{1}(\Omega ))^{N}.  \label{w319}
\end{equation}%
By (\ref{w291}) we can write
\begin{equation*}
\int_{\Omega }(j_{\lambda }(t,x,\nabla u_{\lambda })dx+j_{\lambda }^{\ast
}(t,x,\beta _{\lambda }(t,x,\nabla u_{\lambda }))-\beta _{\lambda
}(t,x,\nabla u_{\lambda })\cdot \nabla u_{\lambda })dx=0,
\end{equation*}%
whence, by (\ref{w307}) and (\ref{w312}), we get that
\begin{eqnarray}
&&\int_{\Omega }(j(t,x,(1+\lambda \beta (t,x,\cdot ))^{-1}\nabla u_{\lambda
})+j_{\lambda }^{\ast }(t,x,\beta _{\lambda }(t,x,\nabla u_{\lambda })))dx
\label{w320} \\
&\leq &\int_{\Omega }\beta _{\lambda }(t,x,\nabla u_{\lambda })\cdot \nabla
u_{\lambda }dx  \notag \\
&\leq &\frac{1}{h}\left\{ \int_{\Omega }w_{1}u_{\lambda }dx+\int_{\Gamma
}w_{2}u_{\lambda }d\sigma -\int_{\Omega }u_{\lambda }^{2}dx-\int_{\Gamma
}u_{\lambda }^{2}d\sigma -\lambda \int_{\Omega }\left\vert \nabla u_{\lambda
}\right\vert _{N}^{2}dx\right\} \leq C.  \notag
\end{eqnarray}%
Passing to the limit in (\ref{w322}), recalling that $(1+\lambda \beta
(t,x,\cdot ))^{-1}\nabla u_{\lambda }\rightarrow \nabla u^{\ast }$ weakly,
we obtain on the basis of the weak lower semicontinuity of the convex
integrand that
\begin{equation}
j(t,\cdot ,\nabla u^{\ast })\in L^{1}(\Omega ).  \label{w323}
\end{equation}%
Also, by (\ref{w320}) it follows that
\begin{equation}
\int_{\Omega }j_{\lambda }^{\ast }(t,x,\beta _{\lambda }(t,x,\nabla
u_{\lambda }))dx\leq C,\mbox{ }\forall \lambda >0.  \label{w324}
\end{equation}%
This yields (by the definition (\ref{w307}) for $j_{\lambda }^{\ast })$%
\begin{multline}
\int_{\Omega }\frac{1}{2\lambda }\left\vert (1+\lambda \beta ^{-1}(t,x,\cdot
))^{-1}\beta _{\lambda }(t,x,\nabla u_{\lambda })-\beta _{\lambda
}(t,x,\nabla u_{\lambda })\right\vert _{N}^{2}  \label{w324-0} \\
+\int_{\Omega }j^{\ast }(t,x,(1+\lambda \beta ^{-1}(t,x,\cdot ))^{-1}\beta
_{\lambda }(t,x,\nabla u_{\lambda }))dx  \notag \\
\leq \int_{\Omega }j_{\lambda }^{\ast }(t,x,\beta _{\lambda }(t,x,\nabla
u_{\lambda }))dx\leq C.  \notag
\end{multline}

Arguing as above, on the basis of (\ref{w294}) and the Dunford-Pettis
theorem we deduce that the sequence $\left\{ (1+\lambda \beta
^{-1}(t,x,\cdot ))^{-1}\beta _{\lambda }(t,x,\nabla u_{\lambda })\right\}
_{\lambda >0}$ is weakly compact in $(L^{1}(\Omega ))^{N},$ and so, on a
subsequence, as $\lambda \rightarrow 0,$ we get
\begin{equation}
(1+\lambda \beta ^{-1}(t,x,\cdot ))^{-1}\beta _{\lambda }(t,x,\nabla
u_{\lambda })\rightarrow \eta \mbox{ weakly in }(L^{1}(\Omega ))^{N},\mbox{ }
\label{w324-1}
\end{equation}%
\begin{equation*}
(1+\lambda \beta ^{-1}(t,x,\cdot ))^{-1}\beta _{\lambda }(t,x,\nabla
u_{\lambda })-\beta _{\lambda }(t,x,\nabla u_{\lambda })\rightarrow 0\mbox{
strongly in }(L^{2}(\Omega ))^{N},
\end{equation*}%
which implies
\begin{equation}
\beta _{\lambda }(t,x,\nabla u_{\lambda })\rightarrow \eta \mbox{ weakly in }%
(L^{1}(\Omega ))^{N}.  \label{w325}
\end{equation}%
Then, by (\ref{w324}) and the weak lower semicontinuity of the convex
integrand we infer that
\begin{equation}
j^{\ast }(t,\cdot ,\eta )\in L^{1}(\Omega ).  \label{w325-0}
\end{equation}%
Now, we pass to the limit in (\ref{w320}), taking into account (\ref{w314})
and (\ref{w324-1}) and we get
\begin{eqnarray}
&&\limsup\limits_{\lambda \rightarrow 0}\int_{\Omega }\beta _{\lambda
}(t,x,\nabla u_{\lambda })\cdot \nabla u_{\lambda }dx  \label{w326} \\
&\leq &\frac{1}{h}\left\{ \int_{\Omega }w_{1}u^{\ast }dx+\int_{\Gamma
}w_{2}u^{\ast }d\sigma -\int_{\Omega }(u^{\ast })^{2}dx-\int_{\Gamma
}(u^{\ast })^{2}d\sigma \right\} .  \notag
\end{eqnarray}

Then, letting $\lambda \rightarrow 0$ in (\ref{w310}) and recalling (\ref%
{w326-0}) we obtain
\begin{equation*}
\int_{\Omega }u^{\ast }\psi dx+h\int_{\Omega }\eta \cdot \nabla \psi
dx+\int_{\Gamma }u^{\ast }\psi d\sigma =\int_{\Omega }w_{1}\psi
dx+\int_{\Gamma }w_{2}\psi d\sigma ,\mbox{ }\forall \psi \in C^{1}(\overline{%
\Omega }).
\end{equation*}%
This is extended by density for all $\psi \in L^{2}(\Omega )\cap
W^{1,1}(\Omega ),$ $\gamma (\psi )\in L^{2}(\Gamma ),$ and in particular for
$\psi =u^{\ast }.$ We obtain
\begin{equation}
\int_{\Omega }\eta \cdot \nabla u^{\ast }dx=\frac{1}{h}\left\{ \int_{\Omega
}w_{1}u^{\ast }dx+\int_{\Gamma }w_{2}u^{\ast }d\sigma -\int_{\Omega
}(u^{\ast })^{2}dx-\int_{\Gamma }(u^{\ast })^{2}d\sigma \right\} .
\label{w327}
\end{equation}%
By (\ref{w326}) and (\ref{w327}) we finally obtain that
\begin{equation}
\limsup\limits_{\lambda \rightarrow 0}\int_{\Omega }\beta _{\lambda
}(t,x,\nabla u_{\lambda })\cdot \nabla u_{\lambda }dx\leq \int_{\Omega }\eta
\cdot \nabla u^{\ast }dx.  \label{w327-1}
\end{equation}%
Since $\nabla u_{\lambda }\rightarrow \nabla u^{\ast }$ weakly in $%
(L^{1}(\Omega ))^{N},$ $\beta _{\lambda }(t,x,\nabla u_{\lambda
})\rightarrow \eta $ weakly in $(L^{1}(\Omega ))^{N},$ we deduce that
\begin{equation*}
\eta (x)\in \beta (t,x,\nabla u^{\ast }(x))\mbox{ a.e. }x\in \overline{%
\Omega }.
\end{equation*}

As a matter of fact, to get the latter we note that by (\ref{w327-1}) we
have
\begin{equation*}
\int_{\Omega }(j_{\lambda }(t,x,\nabla u_{\lambda })-j_{\lambda }(t,x,\theta
))dx\leq \int_{\Omega }\eta \cdot (\nabla u^{\ast }-\theta )dx,\mbox{ }%
\forall \theta \in (L^{1}(\Omega ))^{N},
\end{equation*}%
and letting $\lambda \rightarrow 0$ we get by (\ref{w307}), (\ref{w319}) and
(\ref{w307-0})
\begin{equation*}
\int_{\Omega }(j(t,x,\nabla u^{\ast })-j(t,x,\theta ))dx\leq \int_{\Omega
}\eta \cdot (\nabla u^{\ast }-\theta )dx,\mbox{ }\forall \theta \in
(L^{1}(\Omega ))^{N}.
\end{equation*}%
Since $\theta $ is arbitrary we get that $\eta (x)\in \partial j(t,x,\nabla
u^{\ast }(x)),$ as desired.

Passing to the limit in (\ref{w310}) we get
\begin{equation*}
\int_{\Omega }(u^{\ast }\psi +h\eta \cdot \nabla \psi )dx+\int_{\Gamma
}u^{\ast }\psi d\sigma =b(\psi ),\mbox{ }\forall \psi \in W^{1,2}(\Omega ),
\end{equation*}%
where $b$ is defined by (\ref{w20-00}) for all $u\in U_{1}.$ By density this
extends to all $\psi \in L^{2}(\Omega )\cap W^{1,1}(\Omega ),$ with $\gamma
(\psi )\in L^{2}(\Gamma ).$ Hence $u^{\ast }$ is the weak solution to (\ref%
{w200}). Uniqueness of $u^{\ast },$ as weak solution, is immediate.

Moreover, it also follows that $\nabla \cdot \eta \in L^{2}(\Omega ).$\hfill
$\square $

\medskip

\noindent \textbf{Definition 4.3. }Let
\begin{equation}
y_{0}\in U_{1},\mbox{ }f\in L^{2}(Q),\mbox{ }g\in L^{2}(\Sigma ).
\label{w296-00}
\end{equation}%
We call a \textit{weak solution} to problem (\ref{w1})-(\ref{w3}) a function
$y\in L^{2}(Q),$ such that
\begin{eqnarray}
y &\in &L^{1}(0,T;W^{1,1}(\Omega )),\mbox{ }\gamma (y)\in L^{2}(\Sigma ),%
\mbox{ }j(\cdot ,\cdot ,\nabla y)\in L^{1}(Q),  \label{w296-1} \\
\mbox{and there exists }\eta &\in &(L^{1}(Q))^{N},\mbox{ }\eta (t,x)\in
\beta (t,x,\nabla y(t,x)),\mbox{ a.e. }(t,x)\in Q,  \notag \\
j^{\ast }(\cdot ,\cdot ,\eta ) &\in &L^{1}(Q),\mbox{ such that}  \notag
\end{eqnarray}%
\begin{eqnarray}
&&-\int_{Q}y\phi _{t}dxdt+\int_{Q}\eta \cdot \nabla \phi dxdt-\int_{\Sigma
}y\phi _{t}d\sigma dt  \label{w297} \\
&=&\int_{Q}f\phi dxdt+\int_{\Omega }y_{0}\phi (0)dx+\int_{\Gamma }y_{0}\phi
(0)d\sigma +\int_{\Sigma }g\phi d\sigma dt,  \notag
\end{eqnarray}%
for all $\phi \in W^{1,2}([0,T];L^{2}(\Omega ))\cap L^{1}(0,T;W^{1,1}(\Omega
)),$ $\gamma (\phi )\in W^{1,2}([0,T];L^{2}(\Gamma )),$ with $\phi (T)=0.$

\medskip

By Lemma 4.1 it is clear that the second term on the left-hand side of (\ref%
{w297}) makes sense.

\medskip

\noindent \textbf{Theorem 4.4. }\textit{Let us assume }(\ref{w296-00})%
\textit{\ and }$j(0,\cdot ,\nabla y_{0})\in L^{1}(\Omega )$\textit{. Then,
under hypotheses }(H$_{1}$), (H$_{3})$-(H$_{5}),$ \textit{problem}\textbf{%
\textit{\ }}(\ref{w1})-(\ref{w3})\textbf{\ }\textit{has at least one weak
solution. Moreover, }$y$\textit{\ is a strong solution to }(\ref{w1})-(\ref%
{w3})\textit{, that is, it satisfies }(\ref{w102})-(\ref{w104}).\textit{\
Finally, }$y$ \textit{is given by }
\begin{equation}
y=\lim_{h\rightarrow 0}y^{h}\mbox{\textit{strongly in} }L^{1}(Q),
\label{w330-0}
\end{equation}%
\textit{where }$y^{h}$\textit{\ is defined by }(\ref{w31}). \textit{The
solution} \textit{is unique in the class of functions satisfying }(\ref%
{w296-1}), (\ref{w102})-(\ref{w104})\textit{.}

\medskip

\noindent \textbf{Proof. }Let us consider the time discretized system (\ref%
{w12})-(\ref{w14}) whose weak solution is defined as in Definition 2.1, by
replacing $U$ by $U_{1}.$ We claim that system (\ref{w12})-(\ref{w13}) has a
unique weak solution which satisfies
\begin{eqnarray}
&&\left\Vert y_{m}^{h}\right\Vert _{L^{2}(\Omega )}+\left\Vert \gamma
(y_{m}^{h})\right\Vert _{L^{2}(\Gamma )}  \label{w331} \\
&&+h\sum_{i=0}^{m-1}\int_{\Omega }j(t_{i+1},x,\nabla
y_{i+1}^{h})dx+h\sum_{i=0}^{m-1}\left\Vert \nabla y_{i+1}^{h}\right\Vert
_{L^{1}(\Omega )}\leq C,\mbox{ \ }m=1,...,n,  \notag
\end{eqnarray}%
\begin{equation}
h\sum_{i=0}^{m-1}\left\Vert \frac{y_{i+1}^{h}-y_{i}^{h}}{h}\right\Vert
_{L^{2}(\Omega )}^{2}+h\sum_{i=0}^{m-1}\left\Vert \frac{\gamma
(y_{i+1}^{h})-\gamma (y_{i}^{h})}{h}\right\Vert _{L^{2}(\Gamma )}^{2}\leq C,%
\mbox{ \ }m=1,...,n,  \label{w332}
\end{equation}%
where $C$\ is a positive constant, independent of $h$. The proof follows as
in Proposition 2.2 (see (\ref{w24-0})), using the hypothesis corresponding
to the weakly coercive case.

Next, we define $y^{h}$ by (\ref{w31}) and on the basis of (\ref{w331}) and (%
\ref{w332}) we write
\begin{equation}
\left\Vert y^{h}(t)\right\Vert _{L^{2}(\Omega )}+\left\Vert \gamma
(y^{h}(t))\right\Vert _{L^{2}(\Gamma )}+\int_{Q}j(t,x,\nabla
y^{h}(t))dxdt\leq C,\mbox{ for }t\in \lbrack 0,T],  \label{w333}
\end{equation}%
\begin{equation}
\int_{0}^{T}\left\Vert \frac{y^{h}(t+h)-y^{h}(t)}{h}\right\Vert
_{L^{2}(\Omega )}^{2}dt+\int_{0}^{T}\left\Vert \frac{\gamma
(y^{h}(t+h))-\gamma (y^{h}(t))}{h}\right\Vert _{L^{2}(\Gamma )}^{2}dt\leq C.
\label{w335}
\end{equation}

By these estimates and the Dunford-Pettis compactness theorem in $L^{1}(Q),$
we can select a subsequence such that, as $h\rightarrow 0,$
\begin{equation*}
y^{h}\rightarrow y\mbox{ weak-star in }L^{\infty }(0,T;L^{2}(\Omega )),\mbox{
}
\end{equation*}%
\begin{equation*}
y^{h}\rightarrow y\mbox{ weak-star in }L^{\infty }(0,T;L^{2}(\Gamma )),\mbox{
}
\end{equation*}%
\begin{equation*}
\nabla y^{h}\rightarrow \nabla y\mbox{ weakly in }(L^{1}(Q))^{N},\mbox{ }
\end{equation*}%
\begin{equation*}
\frac{y^{h}(t+h)-y^{h}(t)}{h}\rightarrow \frac{dy}{dt}\mbox{ weakly in }%
L^{2}(Q),\mbox{ }
\end{equation*}%
\begin{equation*}
\frac{\gamma (y^{h}(t+h))-\gamma (y^{h}(t))}{h}\rightarrow \frac{dy}{dt}%
\mbox{ weakly in }L^{2}(\Sigma ).
\end{equation*}

We denote $X=W^{1,r}(\Omega )$ with $r>N.$ Then $W^{1,r}(\Omega )$ is
compact in $L^{\infty }(\Omega )$ and $L^{1}(\Omega )$ is compact in $%
X^{\prime }=(W^{1,r}(\Omega ))^{\prime }.$ We have that $y^{h}\in
BV([0,T];L^{2}(\Omega ))$ which implies that $y^{h}\in BV([0,T];X^{\prime
}). $ By the Helly theorem it follows that
\begin{equation*}
y^{h}(t)\rightarrow y(t)\mbox{ strongly in }X^{\prime }\mbox{, uniformly
with }t\in \lbrack 0,T].
\end{equation*}%
Using again Lemma 5.1 in \cite{vbp-2012} and taking into account that $%
W^{1,1}(\Omega )\subset L^{1}(\Omega )\subset X^{\prime },$ we deduce that
\begin{equation}
y^{h}\rightarrow y\mbox{ strongly in }L^{1}(Q),\mbox{ as }h\rightarrow 0.
\label{w336}
\end{equation}%
The remainder of the proof follows as in Theorem 3.2 and so it will be
omitted.\hfill $\square $

\medskip

\noindent \textbf{Remark 4.5. }The singular case $\beta (t,x,r)\equiv \rho $
sgn $r$ (which is relevant in the study of diffusion systems with singular
energy) is ruled out by the present approach, but, as seen later, the
corresponding problem (\ref{w1})-(\ref{w3}) is well posed, however, in the
space of functions with bounded variation on $\Omega .$

\section{The semigroup approach}

\setcounter{equation}{0}

Everywhere in the following we assume that either hypotheses $($H$_{1}),$ (H$%
_{2}),$ (H$_{5}),$ or (H$_{1}$), (H$_{3}$)-(H$_{5}$) are satisfied. In other
words, we are in one of the cases considered before: strongly coercive or
weakly coercive. Moreover, we assume that $\beta $ is independent of $t,$ $%
\beta \equiv \beta (x,r).$ It shall turn out that in this time-invariant
case Theorems 3.2 and 4.4 can be derived by the nonlinear contraction
semigroup theory which leads to sharper regularity results for the solution $%
y.$

Namely, on the space $X=L^{2}(\Omega )\times L^{2}(\Gamma ),$ endowed with
the standard Hilbertian structure, we consider the operator $\mathcal{A}:D(%
\mathcal{A})\subset X\rightarrow X,$ defined by
\begin{equation}
\mathcal{A}\left(
\begin{array}{c}
u \\
z%
\end{array}%
\right) =\left(
\begin{array}{c}
-\nabla \cdot \beta (x,\nabla u) \\
\beta (x,\nabla u)\cdot \nu%
\end{array}%
\right) ,\mbox{ \ }\forall \left(
\begin{array}{c}
u \\
z%
\end{array}%
\right) \in D(\mathcal{A}),  \label{w418}
\end{equation}%
\begin{equation}
D(\mathcal{A)=}\left\{
\begin{array}{l}
\left(
\begin{array}{c}
u \\
z%
\end{array}%
\right) \in X;\mbox{ }u\in \widetilde{U},\mbox{ }z=\gamma (u),\mbox{ }%
\exists \eta (x)\in \beta (x,\nabla u(x))\mbox{ a.e. }x\in \Omega , \\
\eta \in (L^{1}(\Omega ))^{N},\mbox{ }\nabla \cdot \eta \in L^{2}(\Omega ),%
\mbox{ }\eta \cdot \nu \in L^{2}(\Gamma )%
\end{array}%
\right\} .  \label{w419}
\end{equation}%
Here, $\widetilde{U}=U$ in the strongly coercive case, that is under
hypothesis (H$_{2})$ and $\widetilde{U}=U_{1}$ in the weakly coercive case
under the hypotheses (H$_{3})$-(H$_{4}).$

In (\ref{w418}), by $\beta (x,\nabla u)$ we mean, as usually, any measurable
section $\eta $ of $\beta (x,\nabla u)$ satisfying (\ref{w419})$.$ Then, the
system
\begin{equation}
y_{t}-\nabla \cdot \beta (x,\nabla y)\ni f\mbox{ \ in }Q,  \label{w1-1}
\end{equation}%
\begin{equation}
\beta (x,\nabla y)\cdot \nu +y_{t}\ni g\mbox{\ on }\Sigma ,  \label{w2-1}
\end{equation}%
\begin{equation}
y(0)=y_{0}\mbox{ \ in }\Omega ,  \label{w3-1}
\end{equation}%
can be written as
\begin{eqnarray}
\frac{d}{dt}\left(
\begin{array}{c}
y(t) \\
z(t)%
\end{array}%
\right) +\mathcal{A}\left(
\begin{array}{c}
y(t) \\
z(t)%
\end{array}%
\right) &\ni &\left(
\begin{array}{c}
f(t) \\
g(t)%
\end{array}%
\right) ,\mbox{ \ a.e. }t\in (0,T),  \label{w420} \\
\left(
\begin{array}{c}
y \\
z%
\end{array}%
\right) (0) &=&\left(
\begin{array}{c}
y_{0} \\
z_{0}%
\end{array}%
\right) .  \notag
\end{eqnarray}

\medskip

\noindent \textbf{Lemma 5.1. }\textit{The operator }$\mathcal{A}$\textit{\
is maximal monotone in }$X.$

\medskip

\noindent \textbf{Proof. }It is easily seen that $\mathcal{A}$ is monotone,
that is,
\begin{equation*}
\left( \mathcal{A}\left(
\begin{array}{c}
u \\
z%
\end{array}%
\right) -\mathcal{A}\left(
\begin{array}{c}
\overline{u} \\
\overline{z}%
\end{array}%
\right) ,\left(
\begin{array}{c}
u-\overline{u} \\
z-\overline{z}%
\end{array}%
\right) \right) _{X}\geq 0,\mbox{ \ }\forall \left(
\begin{array}{c}
u \\
z%
\end{array}%
\right) ,\left(
\begin{array}{c}
\overline{u} \\
\overline{z}%
\end{array}%
\right) \in D(\mathcal{A}).
\end{equation*}%
In fact, this follows by the Gauss-Ostrogradski formula
\begin{equation*}
-\int_{\Omega }v\nabla \cdot \eta dx=-\int_{\Gamma }\gamma (v)\eta \cdot \nu
d\sigma +\int_{\Omega }\eta \cdot \nabla vdx,\mbox{ }\eta (x)\in \beta
(x,\nabla u(x)),\mbox{ a.e. }x\in \Omega ,
\end{equation*}%
for all $u,v\in U$, in the strongly coercive case, or $u,v\in U_{1}$ in the
weakly coercive case.

On the other hand, the range $R(I+\mathcal{A})$ is all of $X,$ for all $%
\lambda >0.$ Indeed, equation
\begin{equation*}
(I+\mathcal{A})\left(
\begin{array}{c}
u \\
z%
\end{array}%
\right) =\left(
\begin{array}{c}
w_{1} \\
w_{2}%
\end{array}%
\right) \mbox{ for all }\left(
\begin{array}{c}
w_{1} \\
w_{2}%
\end{array}%
\right) \in X
\end{equation*}
reduces to equation (\ref{w20-0}) (or (\ref{w200})), for which existence has
been previously proved.\hfill $\square $

\medskip

Then, by the standard existence theorem for the Cauchy problem associated
with nonlinear maximal monotone operators (see, e.g., \cite{vb-springer-2010}%
, p. 151), for $\left(
\begin{array}{c}
y_{0} \\
z_{0}%
\end{array}%
\right) \in D(\mathcal{A})$ and $f\in W^{1,1}([0,T];L^{2}(\Omega )),$ $g\in
W^{1,1}([0,T];L^{2}(\Gamma )),$ there is a unique function $\left(
\begin{array}{c}
y \\
z%
\end{array}%
\right) \in W^{1,\infty }([0,T];L^{2}(\Omega ))\times W^{1,\infty
}([0,T];L^{2}(\Gamma ))$ which satisfies (\ref{w420}) a.e., and also in the
following stronger sense,%
\begin{equation*}
\frac{d^{+}}{dt}\left(
\begin{array}{c}
y(t) \\
z(t)%
\end{array}%
\right) +\left( \mathcal{A}\left(
\begin{array}{c}
y(t) \\
z(t)%
\end{array}%
\right) -\left(
\begin{array}{c}
f(t) \\
g(t)%
\end{array}%
\right) \right) ^{\circ }=0,\mbox{ \ }\forall t\in \lbrack 0,T),
\end{equation*}%
where, for each closed convex set $C,$ $C^{\circ }$ stands for the minimal
section of $C$. Moreover, if $f=0,$ $g=0,$ we have the exponential formula
\begin{equation*}
\left(
\begin{array}{c}
y(t) \\
z(t)%
\end{array}%
\right) =\lim_{n\rightarrow \infty }\left( I+\frac{t}{n}\mathcal{A}\right)
^{-n}\left(
\begin{array}{c}
y_{0} \\
z_{0}%
\end{array}%
\right) \mbox{ in }L^{2}(\Omega )\times L^{2}(\Gamma ),
\end{equation*}%
uniformly in $t$ on compact intervals.

We are led therefore to the following sharper versions of Theorems 3.2 and
4.4.

\medskip

\noindent \textbf{Theorem 5.2. }\textit{Let }$y_{0}\in U$\textit{\
(respectively }$U_{1}$ \textit{in the weakly coercive case}$)$\textit{\ be
such that }$\nabla \cdot \beta (\cdot ,\nabla y_{0})\in L^{2}(\Omega ),$%
\textit{\ }$\beta (\cdot ,\nabla y_{0})\cdot \nu \in L^{2}(\Gamma ),$\textit{%
\ and let }$f\in W^{1,1}([0,T];L^{2}(\Omega )),$\textit{\ }$g\in
W^{1,1}([0,T];L^{2}(\Gamma )).$\textit{\ Then, there is a unique }$y\in
W^{1,\infty }([0,T];L^{2}(\Omega ))$\textit{\ with }$\gamma (y)\in
W^{1,\infty }([0,T];L^{2}(\Gamma ))$\textit{\ which satisfies }%
\begin{equation*}
\frac{d^{+}}{dt}y(t,x)-(\nabla \cdot \beta (x,\nabla y(t,x))-f(t,x))^{\circ
}=0\mbox{ in }[0,T)\times \Omega ,
\end{equation*}%
\begin{equation*}
\frac{d^{+}}{dt}y(t,x)+(\beta (x,\nabla y(t,x))\cdot \nu -g(t,x))^{\circ }=0%
\mbox{ on }[0,T)\times \Gamma ,
\end{equation*}%
\begin{equation*}
y(0,x)=y_{0}\mbox{ in }\Omega .
\end{equation*}%
\textit{Moreover, }(\ref{w105})\textit{\ holds.}

\medskip

Condition $\nabla \cdot \beta (\cdot ,\nabla y_{0})\in L^{2}(\Omega )$
means, of course, that there is $\eta _{0}$ measurable, such that $\eta
_{0}(x)\in \beta (x,\nabla y_{0}(x)),$ a.e. $x\in \Omega ,$ $\nabla \cdot
\eta _{0}\in L^{2}(\Omega ).$

\noindent We also note that the operator $\mathcal{A}$ is the
subdifferential of the function $\Phi :X\rightarrow \mathbb{R}\cup \{+\infty
\},$%
\begin{equation*}
\Phi \left(
\begin{array}{c}
u \\
z%
\end{array}%
\right) =\left\{
\begin{array}{l}
\int_{\Omega }j(x,\nabla u(x))dx,\mbox{ }z=\gamma (u)\in L^{2}(\Gamma )\mbox{%
\ \ \ if }u\in W^{1,1}(\Omega ),\mbox{ }j(\cdot ,\nabla u)\in L^{1}(\Omega )
\\
\mbox{ \ \ } \\
+\infty ,\mbox{ \ \ \ \ \ \ \ \ \ \ \ \ \ \ \ \ \ \ \ \ \ \ \ \ \ \ \ \ \ \
\ \ \ \ \ \ \ \ \ \ \ \ otherwise.}%
\end{array}%
\right.
\end{equation*}%
This is the energy functional associated with system (\ref{w1-1})-(\ref{w3-1}%
).

Indeed, for $\left(
\begin{array}{c}
u \\
z%
\end{array}%
\right) \in D(\mathcal{A}),\left(
\begin{array}{c}
\overline{u} \\
\overline{z}%
\end{array}%
\right) \in D(\Phi )$ we have
\begin{eqnarray*}
&&\Phi \left(
\begin{array}{c}
u \\
z%
\end{array}%
\right) -\Phi \left(
\begin{array}{c}
\overline{u} \\
\overline{z}%
\end{array}%
\right) =\int_{\Omega }(j(x,\nabla u)-j(x,\nabla \overline{u}))dx \\
&\leq &\int_{\Omega }\eta \cdot (\nabla u-\nabla \overline{u}%
)dx=-\int_{\Omega }(u-\overline{u})\nabla \cdot \eta dx+\int_{\Gamma }(\eta
\cdot \nu )(z-\overline{z})d\sigma ,
\end{eqnarray*}%
for $\eta \in (L^{1}(\Omega ))^{N},$ $\eta (x)\in \beta (x,\nabla u(x))$
a.e. $x\in \Omega .$

Here we have used the Gauss-Ostrogradski formula%
\begin{equation*}
\int_{\Omega }v\cdot \nabla udx=-\int_{\Omega }u\nabla \cdot
vdx+\int_{\Gamma }\gamma (u)v\cdot \nu d\sigma ,
\end{equation*}%
which is valid for all $u\in W^{1,1}(\Omega )\cap L^{2}(\Omega )$ and $v\in
(L^{1}(\Omega ))^{N}$ such that $\gamma (u)(v\cdot \nu )\in L^{1}(\Gamma )$
and $u\nabla \cdot v\in L^{1}(\Omega ).$ In virtue of Lemma 4.1, $v=\eta $
satisfies this condition.

This implies that $\mathcal{A}\subset \partial \Phi $ and since $\mathcal{A}$
is maximal monotone we infer that $A=\partial \Phi ,$ as claimed.

Then, by Theorem 4.11 in \cite{vb-springer-2010}, p. 158, it follows that
\textit{for all }$\left(
\begin{array}{c}
y_{0} \\
z_{0}%
\end{array}%
\right) \in D(\Phi ),$\textit{\ }$f\in L^{2}(Q),$\textit{\ }$g\in
L^{2}(\Sigma ),$\textit{\ problem }(\ref{w420})\textit{\ has a unique
solution }$\left(
\begin{array}{c}
y \\
z%
\end{array}%
\right) \in W^{1,2}([0,T];X).$ Hence, under the above assumptions \textit{%
there is a unique solution }$y\in W^{1,2}([0,T];L^{2}(\Omega ))$\textit{\
with }$\gamma (y)\in W^{1,2}([0,T];L^{2}(\Gamma ))$\textit{\ to} (\ref{w1-1}%
)-(\ref{w3-1}).

By Theorem 4.13 in \cite{vb-springer-2010}, p. 164, we have also in this
case the following asymptotic result for the solution $y$ to (\ref{w1-1})-(%
\ref{w3-1}).

\medskip

\noindent \textbf{Theorem 5.3. }\textit{Let }$y_{0}\in W^{1,p}(\Omega ),$%
\textit{\ }$2\leq p<\infty $ and $f(t)\equiv f\in L^{2}(\Omega ),$ $%
g(t)\equiv g\in L^{2}(\Gamma ).$ \textit{Assume that the set of equilibrium
states for }(\ref{w1}),\textit{\ }%
\begin{multline*}
\mathcal{K}=\{y\in U;\nabla \cdot \beta (\cdot ,\nabla y)\in L^{2}(\Omega ),%
\mbox{ }\beta (\cdot ,\nabla y)\cdot \nu \in L^{2}(\Gamma ),\mbox{ } \\
\nabla \cdot \beta (x,\nabla y)=f\mbox{ in }\Omega ,\mbox{ }\beta (x,\nabla
y)\cdot \nu =g\mbox{ on }\Gamma \}
\end{multline*}%
\textit{is non empty. Then, for }$t\rightarrow \infty ,$\textit{\ we have}
\begin{eqnarray}
y(t) &\rightarrow &y_{\infty }\mbox{ \textit{weakly in} }L^{2}(\Omega ),
\label{w420-0} \\
\gamma (y(t)) &\rightarrow &\gamma (y_{\infty })\mbox{ \textit{weakly in} }%
L^{2}(\Gamma ),  \notag
\end{eqnarray}%
where $y_{\infty }\in \mathcal{K}.$

\medskip

Taking into account that, as easily seen by (\ref{w420}) we have
\begin{equation*}
\Phi \left(
\begin{array}{c}
y(t) \\
z(t)%
\end{array}%
\right) \leq \Phi \left(
\begin{array}{c}
y_{0} \\
\gamma (y_{0})%
\end{array}%
\right) ,\mbox{ }\forall t\geq 0,
\end{equation*}%
it follows by (\ref{w420-0}) and by the compactness of $W^{1,p}(\Omega )$ in
$L^{q}(\Omega ),$ that for $t\rightarrow \infty ,$%
\begin{equation*}
y(t)\rightarrow y_{\infty }\mbox{ strongly in }L^{q}(\Omega ),
\end{equation*}%
where $1\leq q<\frac{Np}{N-p}$ if $N>p,$ $q=p$ if $N\leq p$ in the strongly
coercive case, and $q=1$ in the weakly coercive case.

In other words, the solution $y$ is strongly convergent to an equilibrium
solution $y_{\infty }$ to system (\ref{w1-1})-(\ref{w3-1}).

One of the main advantages of the semigroup approach is its flexibility to
incorporate other nonlinear terms in the basic equations (\ref{w1-1})-(\ref%
{w2-1}). We shall consider two such extensions. The first is the problem
studied in \cite{warma-2012}, already mentioned in Introduction,\textit{\ }%
\begin{equation}
\frac{\partial y}{\partial t}-\nabla \cdot \beta (x,\nabla y)+a_{1}(x,y)\ni
f(t,x)\mbox{ in }Q,  \label{w421}
\end{equation}%
\begin{equation}
\frac{\partial y}{\partial t}-\nabla \cdot (\left\vert \nabla _{\Gamma
}y\right\vert _{N-1}^{p-2}\nabla _{\Gamma }y)+\beta (x,\nabla y)\cdot \nu
+a_{2}(x,y)\ni g(t,x)\mbox{ on }\Sigma ,  \label{w422}
\end{equation}%
\begin{equation}
y(0,x)=y_{0}\mbox{ in }\Omega ,  \label{w423}
\end{equation}%
where $\beta $ satisfies assumption (\ref{w5})\textit{, }$a_{i}:\overline{%
\Omega }\times \mathbb{R}\rightarrow \mathbb{R},$ $i=1,2$ are continuous and
$p\geq 2.$ Here, $\nabla _{\Gamma }y$ is the Riemannian gradient of $y,$
that is $\nabla _{\Gamma }y=(\partial _{\tau _{1}}y,...,\partial _{\tau
_{N-1}}y),$ where $\partial _{\tau _{i}}y$ is the directional derivative of $%
y$ along the tangential directions $\tau _{i}$ at each point on $\Gamma $
(see \cite{warma-2012}) and%
\begin{equation*}
-\int_{\Gamma }v\nabla \cdot (\left\vert \nabla _{\Gamma }u\right\vert
^{p-2}\nabla _{\Gamma }u)d\sigma =\int_{\Gamma }\left\vert \nabla _{\Gamma
}u\right\vert ^{p-2}\nabla _{\Gamma }y\cdot \nabla _{\Gamma }vd\sigma .
\end{equation*}

Problem (\ref{w421})-(\ref{w423}) can be written as (\ref{w420}), where%
\begin{equation}
\mathcal{A}\left(
\begin{array}{c}
u \\
z%
\end{array}%
\right) =\left(
\begin{array}{c}
-\nabla \cdot \beta (x,\nabla u)+a_{1}(x,u) \\
-\nabla \cdot (\left\vert \nabla _{\Gamma }u\right\vert _{N-1}^{p-2}\nabla
_{\Gamma }u)+\beta (x,\nabla u)\cdot \nu +a_{2}(x,u)%
\end{array}%
\right) ,\mbox{ \ }  \label{w441}
\end{equation}%
for all $\left(
\begin{array}{c}
u \\
z%
\end{array}%
\right) \in D(\mathcal{A}),$ where
\begin{equation}
D(\mathcal{A)=}\left\{
\begin{array}{l}
\left(
\begin{array}{c}
u \\
z%
\end{array}%
\right) \in X;\mbox{ }u\in U,\mbox{ }z=\gamma (u),\mbox{ }\exists \eta
(x)\in \beta (x,\nabla u(x))\mbox{ a.e. }x\in \Omega ,\mbox{ such that} \\
\eta \in (L^{1}(\Omega ))^{N},\mbox{ }\nabla \cdot \eta +a_{1}(\cdot ,u)\in
L^{2}(\Omega ),\mbox{ }\left\vert \nabla _{\Gamma }u\right\vert _{N-1}\in
L^{p}(\Gamma ),\mbox{ } \\
-\nabla \cdot (\left\vert \nabla _{\Gamma }u\right\vert _{N-1}^{p-2}\nabla
_{\Gamma }u)+\eta \cdot \nu +a_{2}(\cdot ,u)\in L^{2}(\Gamma )%
\end{array}%
\right\} .  \label{w442}
\end{equation}

If $y\rightarrow a_{i}(x,y),$ $i=1,2,$ are monotone (or more generally
quasi-monotone, that is, $\lambda y+a_{i}(x,y)$ are monotone for some $%
\lambda >0)$ and $\left\vert a_{i}(x,r)\right\vert \leq C\left\vert
r\right\vert ^{q-1},$ $\forall r\in \mathbb{R},$ where $q$ is as before,
then, arguing as above, it follows that the operator $\mathcal{A}$ is
maximal monotone in $X$ (or quasi $m$-accretive if $a_{i}$ are quasi
monotone), and in fact it is a subdifferential operator.

Then, we get for problem (\ref{w421})-(\ref{w423}) the following existence
result: \textit{let }$y_{0}\in U,$\textit{\ such that }$\nabla _{\Gamma
}y_{0}\in (L^{p}(\Gamma ))^{N-1}$\textit{\ and let }$f\in
L^{2}(0,T;L^{2}(\Omega )),$\textit{\ }$g\in L^{2}(0,T;L^{2}(\Gamma )).$%
\textit{\ Then, there is a unique solution }$y\in W^{1,2}([0,T];L^{2}(\Omega
))$\textit{\ to }(\ref{w421})-(\ref{w423})\textit{, such that }$\gamma
(y)\in W^{1,2}([0,T];L^{2}(\Gamma ))$\textit{\ and }$\nabla _{\Gamma }y\in
(L^{p}(\Sigma ))^{N-1}.$

\paragraph{The obstacle problem}

Consider the following free boundary problem associated with the Wentzell
boundary condition, namely,
\begin{eqnarray}
y_{t}-\nabla \cdot \beta (x,\nabla y) &\geq &f,\mbox{ }y\geq 0,\mbox{ in }Q,
\label{w442-0} \\
y_{t}-\nabla \cdot \beta (x,\nabla y) &\ni &f,\mbox{ in }\{(t,x)\in Q;\mbox{
}y(t,x)>0\},  \notag \\
y_{t}+\beta (x,\nabla y)\cdot \nu &\ni &g,\mbox{ on }\Sigma ,  \notag \\
y(0,x) &=&y_{0},\mbox{ in }\Omega .  \notag
\end{eqnarray}

This problem can be written as%
\begin{eqnarray}
\frac{d}{dt}\left(
\begin{array}{c}
y(t) \\
z(t)%
\end{array}%
\right) +\mathcal{A}\left(
\begin{array}{c}
y(t) \\
z(t)%
\end{array}%
\right) +\mathcal{B}\left(
\begin{array}{c}
y(t) \\
z(t)%
\end{array}%
\right) &\ni &\left(
\begin{array}{c}
f(t) \\
g(t)%
\end{array}%
\right) ,\mbox{ \ a.e. }t\in (0,T),  \label{w442-1} \\
\left(
\begin{array}{c}
y \\
z%
\end{array}%
\right) (0) &=&\left(
\begin{array}{c}
y_{0} \\
\gamma (y_{0})%
\end{array}%
\right) ,  \notag
\end{eqnarray}%
where
\begin{equation*}
\mathcal{B}\left(
\begin{array}{c}
y \\
z%
\end{array}%
\right) =\left(
\begin{array}{c}
a(y) \\
0%
\end{array}%
\right) ,\mbox{ }\forall \left(
\begin{array}{c}
y \\
z%
\end{array}%
\right) \in X=L^{2}(\Omega )\times L^{2}(\Gamma )
\end{equation*}%
and $a:\mathbb{R}\rightarrow \mathbb{R}$ is the multivalued function $a(s)=0$
for $s>0,$ $a(0)=(-\infty ,0],$ $a(s)=\varnothing $ for $s<0.$ Taking into
account that $\left( \mathcal{A}\left(
\begin{array}{c}
y \\
z%
\end{array}%
\right) ,\mathcal{B}\left(
\begin{array}{c}
y \\
z%
\end{array}%
\right) \right) _{X}\geq 0$ for all $\left(
\begin{array}{c}
y \\
z%
\end{array}%
\right) \in D(\mathcal{A})\cap D(\mathcal{B}),$ it follows that $\mathcal{A}+%
\mathcal{B}$ is maximal monotone and therefore $\mathcal{A+B}=\partial \Phi
_{1},$ where
\begin{equation*}
\Phi _{1}\left(
\begin{array}{c}
y \\
z%
\end{array}%
\right) =\left\{
\begin{array}{l}
\Phi \left(
\begin{array}{c}
y \\
z%
\end{array}%
\right) ,\mbox{ for }y\in L^{2}(\Omega ),\mbox{ }y\geq 0\mbox{ a.e. in }%
\Omega , \\
+\infty ,\mbox{ \ \ \ \ \ \ otherwise.}%
\end{array}%
\right.
\end{equation*}%
Then, applying the general existence theory, we infer that \textit{for }$%
y_{0}\in W^{1,1}(\Omega ),$\textit{\ such that }$y_{0}\geq 0$\textit{\ a.e.
in }$\Omega ,$\textit{\ and }$j(\cdot ,\nabla y_{0})\in L^{1}(\Omega ),$%
\textit{\ problem} (\ref{w442-1}) \textit{has a unique solution }$y\in
W^{1,2}([0,T];L^{2}(\Omega )).$

More generally, one might take instead of $a$ a general maximal monotone
graph in $\mathbb{R}\times \mathbb{R}.$ This case is studied in \cite{frggr}.

\paragraph{The total variation Wentzell flow}

Let us consider now the singular case $j(x,r)\equiv \rho \left\vert
r\right\vert ,$ $r\in \mathbb{R}^{N},$ or equivalently
\begin{equation*}
\beta (x,r)=\rho \mbox{ sgn }r=\left\{
\begin{array}{l}
\rho \frac{r}{\left\vert r\right\vert _{N}},\mbox{ \ \ \ \ \ \ \ \ \ \ \ \
if }r\neq 0, \\
\{r;\left\vert r\right\vert _{N}\leq \rho \},\mbox{ \ \ if }r=0.%
\end{array}%
\right.
\end{equation*}%
Then, problem (\ref{w1-1})-(\ref{w3-1}) reduces to%
\begin{equation}
y_{t}-\rho \nabla \cdot \mbox{sgn }(\nabla y)\ni f,\mbox{ in }Q,
\label{w443}
\end{equation}%
\begin{equation}
\rho \mbox{ sgn }(\nabla y)\cdot \nu +y_{t}\ni g,\mbox{\ on }\Sigma ,
\label{w444}
\end{equation}%
\begin{equation}
y(0)=y_{0},\mbox{ in }\Omega .  \label{w445}
\end{equation}%
As mentioned earlier, this problem is not covered by the previous weakly
coercive case and, as a matter of fact, it cannot be treated in the $%
W^{1,1}(\Omega )$ space, but in the space $BV(\Omega )$ of functions with
bounded variation on $\Omega ,$ that is
\begin{equation*}
BV(\Omega )=\left\{ u\in L^{1}(\Omega );\left\Vert Du\right\Vert
=\sup_{\left\Vert \varphi \right\Vert _{L^{\infty }(\Omega ;\mathbb{R}%
^{N})}\leq 1}\left\{ \int_{\Omega }u\nabla \cdot \varphi dx;\varphi \in
C_{0}^{\infty }(\Omega ;\mathbb{R}^{N})\right\} <\infty \right\} .
\end{equation*}

We recall (see e.g., \cite{attouch}) that for each $u\in BV(\Omega )$ there
is the trace $\gamma (u)\in L^{1}(\Gamma ;d\mathcal{H}^{N-1})$, where $d%
\mathcal{H}^{N-1}$ is the Hausdorff measure on $\Gamma ,$ defined by
\begin{equation*}
\int_{\Omega }u\nabla \cdot \psi dx=-\int_{\Omega }\psi d(\nabla
u)+\int_{\Gamma }u\psi \cdot \nu d\mathcal{H}^{N-1},\mbox{ }\forall \psi \in
C^{1}(\mathbb{R}^{N},\mathbb{R}^{N}).
\end{equation*}%
Here, $\nabla u$ (the gradient of $u$ in the sense of distributions$)$ is a
Radon measure on $\Omega .$

Let us define the energy functional $\Phi :L^{2}(\Omega )\times L^{2}(\Gamma
)\rightarrow (-\infty ,+\infty ],$%
\begin{equation*}
\Phi \left(
\begin{array}{c}
u \\
z%
\end{array}%
\right) =\left\{
\begin{array}{l}
\rho \left\Vert Du\right\Vert ,\mbox{ \ if }u\in BV(\Omega )\cap
L^{2}(\Omega ),\mbox{ }z=\gamma (u)\in L^{2}(\Gamma ), \\
+\infty ,\mbox{ \ \ \ \ \ otherwise.}%
\end{array}%
\right.
\end{equation*}%
It is easily seen that $\Phi $ is convex and l.s.c on $X=L^{2}(\Omega
)\times L^{2}(\Gamma ).$ Let $\partial \Phi :X\rightarrow X$ be its
subdifferential. Then, for each $\left(
\begin{array}{c}
y_{0} \\
z_{0}%
\end{array}%
\right) \in D(\Phi )$ and $f\in L^{2}(0,T;L^{2}(\Omega )),$ $g\in
L^{2}(0,T;L^{2}(\Gamma ))$ the problem
\begin{equation}
\frac{d}{dt}\left(
\begin{array}{c}
y \\
z%
\end{array}%
\right) (t)+\left(
\begin{array}{c}
\xi \\
\eta%
\end{array}%
\right) (t)=\left(
\begin{array}{c}
f \\
g%
\end{array}%
\right) (t),\mbox{ a.e. }t\in (0,T),  \label{w446}
\end{equation}%
\begin{equation}
\left(
\begin{array}{c}
\xi (t) \\
\eta (t)%
\end{array}%
\right) \in \partial \Phi \left(
\begin{array}{c}
y(t) \\
z(t)%
\end{array}%
\right) ,\mbox{ a.e. }t\in (0,T),  \label{w447}
\end{equation}%
\begin{equation}
\left(
\begin{array}{c}
y \\
z%
\end{array}%
\right) (0)=\left(
\begin{array}{c}
y_{0} \\
z_{0}%
\end{array}%
\right) ,\mbox{ in }\Omega ,  \label{w448}
\end{equation}%
has a unique solution $\left(
\begin{array}{c}
y \\
z%
\end{array}%
\right) \in W^{1,2}([0,T];X),$ $\left(
\begin{array}{c}
\xi \\
\eta%
\end{array}%
\right) \in L^{2}(0,T;X).$

Taking into account\ (see e.g. \cite{andreu}\ ) that for all $u,v\in
BV(\Omega )\cap L^{2}(\Omega )$ and $\zeta \in (L^{\infty }(\Omega ))^{N},$%
\begin{equation*}
\left\Vert Du\right\Vert \leq \left\Vert Dv\right\Vert -\int_{\Omega
}(u-v)\nabla \cdot \zeta dx-\int_{\Gamma }(\zeta \cdot \nu )(u-v)d\mathcal{H}%
^{N-1},
\end{equation*}%
where $\left\Vert \zeta \right\Vert _{(L^{\infty }(\Omega ))^{N}}\leq 1,$ $%
\nabla \cdot \zeta \in L^{2}(\Omega ),$ we may interpret\ $t\rightarrow y(t)$
as a solution to system (\ref{w443})-(\ref{w445}). This is the total
variation Wentzell flow.

The operator $\partial \Phi $ (and implicitly system (\ref{w446})-(\ref{w448}%
)) is however hardly to be described in explicit terms, so that a better
insight into problem (\ref{w443})-(\ref{w445}) can be gained by taking into
account that the solution to (\ref{w446})-(\ref{w448}) is the limit of the
finite difference scheme provided by the iteration process%
\begin{equation*}
\left(
\begin{array}{c}
y_{i+1} \\
z_{i+1}%
\end{array}%
\right) +h\partial \Phi \left(
\begin{array}{c}
y_{i+1} \\
z_{i+1}%
\end{array}%
\right) \ni \left(
\begin{array}{c}
y_{i} \\
z_{i}%
\end{array}%
\right) ,\mbox{ }i=0,1,
\end{equation*}%
or equivalently%
\begin{equation*}
y_{i+1}=\arg \min_{u}\left\{ \rho h\left\Vert Du\right\Vert +\frac{1}{2}%
\int_{\Omega }\left\vert u-y_{i}\right\vert ^{2}dx+\frac{1}{2}\int_{\Gamma
}\left\vert \gamma (u)-\gamma (y_{i})\right\vert ^{2}d\sigma \right\} .
\end{equation*}

\paragraph{Final remarks}

The previous results naturally extend to the case of nonlinear functions $%
\beta :\overline{\Omega }\times \mathbb{R}^{N}\rightarrow \mathbb{R}^{N},$
which are not of gradient type with respect to $r\in \mathbb{R}^{N}.$
Namely, it suffices to assume that $\beta \equiv \beta (x,r)$ is continuous
on $\overline{\Omega }\times \mathbb{R}^{N},$ monotone with respect to $r,$
that is
\begin{equation}
(\beta (x,r)-\beta (x,\overline{r}))\cdot (r-\overline{r})\geq 0,\mbox{ }%
\forall r,\overline{r}\in \mathbb{R}^{N},  \label{w449}
\end{equation}%
and that it satisfies
\begin{equation}
\beta (x,r)\cdot r\geq \alpha _{1}\left\vert r\right\vert _{N}^{p},\mbox{ }%
\forall r\in \mathbb{R}^{N},  \label{w450}
\end{equation}%
\begin{equation}
\left\vert \beta (x,r)\right\vert _{N}\leq \alpha _{2}\left\vert
r\right\vert _{N}^{p-1}+\alpha _{3},\mbox{ }\forall r\in \mathbb{R}^{N},
\label{w451}
\end{equation}%
where $2\leq p<\infty ,$ $\alpha _{1},\alpha _{2}>0,$ $\alpha _{3}\in
\mathbb{R}.$

Let us consider
\begin{equation*}
\mathcal{U}=\left\{ \left(
\begin{array}{c}
y \\
z%
\end{array}%
\right) \in U\times L^{2}(\Gamma );\mbox{ }\gamma (u)=z\right\}
\end{equation*}%
endowed with the natural norm and denote by $\mathcal{U}^{\prime }$ the dual
space, in the duality induced by the pivot space $X.$ Then, the operator $%
\widetilde{\mathcal{A}}:\mathcal{U}\rightarrow \mathcal{U}^{\prime },$
defined by
\begin{equation*}
\left\langle \widetilde{\mathcal{A}}\left(
\begin{array}{c}
u \\
z%
\end{array}%
\right) ,\left(
\begin{array}{c}
\varphi _{1} \\
\varphi _{2}%
\end{array}%
\right) \right\rangle _{\mathcal{U}^{\prime },\mathcal{U}}=\int_{\Omega
}\beta (x,\nabla u)\cdot \nabla \varphi _{1}dx,\mbox{ }\forall \left(
\begin{array}{c}
\varphi _{1} \\
\varphi _{2}%
\end{array}%
\right) \in \mathcal{U}
\end{equation*}%
is, by the Browder theory (see e.g., \cite{vb-springer-2010}, p. 81),
maximal monotone in $\mathcal{U}\times \mathcal{U}^{\prime }$ and so its
restriction
\begin{equation*}
\mathcal{A}_{X}\left(
\begin{array}{c}
u \\
z%
\end{array}%
\right) =\widetilde{\mathcal{A}}\left(
\begin{array}{c}
u \\
z%
\end{array}%
\right) \cap X
\end{equation*}%
to $X$ is maximal monotone in $X\times X.$ Then, Theorem 5.2 remains true in
the present situation.

\mathstrut

\noindent \textbf{Acknowledgment. }{\small This work was partially supported
by a grant of the Romanian National Authority for Scientific Research, CNCS
--UEFISCDI, project number PN-II-ID-PCE-2011-3-0027 and also by an
Italian-Romanian project within the scientific agreement between CNR Italy
and the Romanian Academy.}

\end{document}